\documentclass[leqno,12pt]{article} 
\usepackage{amsmath, amssymb}
\usepackage{amsthm} 
\usepackage[all]{xy}
\usepackage[ansinew]{inputenc} 
\newcommand{\Norm}{{\mathcal N}}
\newcommand{\F}{{\mathcal F}}

\newcommand{\tl}{\widetilde}
\newcommand{\Ol}{{\mathcal O}}
\newcommand{\f}{\varphi}
\newcommand{\ac}{{\mathcal H}}
\newcommand{\E}{{\mathcal E}}

\newcommand{\V}{{\mathcal V}}
\newcommand{\G}{{\mathbb G}}

\newcommand{\pu}{{\mathbb P^1}}
\newcommand{\proj}{\mathbb P}
\newcommand{\quadr}{\mathbb Q}
\newcommand{\qu}{\mathcal Q}
\newcommand{\scroll}{\mathcal S}
\newcommand{\pf}{{\mathbb P^4}}
\newcommand{\pt}{{\mathbb P^3}}
\newcommand{\pd}{{\mathbb P^2}}
\newcommand{\lra}{{\longrightarrow}}
\DeclareMathOperator{\loc}{\mathrm{Locus}}
\DeclareMathOperator{\cloc}{\mathrm{ChLocus}}
\DeclareMathOperator{\cone}{NE}
\DeclareMathOperator{\Bs}{Bs}
\DeclareMathOperator{\pic}{Pic}
\DeclareMathOperator{\cycl}{N_1}
\DeclareMathOperator{\Bl}{Bl}
\DeclareMathOperator{\Exc}{Exc}
\DeclareMathOperator{\codim}{codim}
\DeclareMathOperator{\rk}{rk}
\DeclareMathOperator{\ord}{ord}
\DeclareMathOperator{\cl}{cl}
\newcommand{\sa}{\hspace{1pt}}
\newcommand{\W}{{\overline{W}}\sa}
\newcommand{\onespan}[1]{\langle #1 \rangle}
\newcommand{\twospan}[2]{\langle #1,#2 \rangle}

\newcommand{\morespan}[2]{\langle #1, \ldots, #2 \rangle}
\newcommand{\oth}{{\overline{\theta}}}

\newcommand{\os}{{\overline{\sigma}}}

\newcommand{\ms}{\medskip}

\newcommand{\ratcurves}{\textrm{Ratcurves}^n(X)}

\newcommand{\shse}[3]{0 ~\lra ~#1~ \lra ~#2~ \lra ~#3~ \lra~ 0}

\renewcommand{\baselinestretch}{1.2} 
\theoremstyle{plain} 
\newtheorem{theorem}{\indent\sc Theorem}[section]
\newtheorem{lemma}[theorem]{\indent\sc Lemma}
\newtheorem{corollary}[theorem]{\indent\sc Corollary}
\newtheorem{proposition}[theorem]{\indent\sc Proposition}

\theoremstyle{definition} 
\newtheorem{definition}[theorem]{\indent\sc Definition}
\newtheorem{remark}[theorem]{\indent\sc Remark}
\newtheorem{example}[theorem]{\indent\sc Example}

\makeatletter
\def\address#1#2{\begingroup
\noindent\parbox[t]{7.8cm}{%
\small{\scshape\ignorespaces#1}\par\vskip1ex
\noindent\small{\itshape E-mail address}%
\/: #2\par\vskip4ex}\hfill%
\endgroup}%
\makeatother
\title{\uppercase{Fano fivefolds of index two with blow-up structure}} 
\author{\textsc{Elena Chierici* and Gianluca Occhetta} 
}
\date{} 
\begin{document}

\maketitle


\footnote{ 
$^{*}$ The first named author has been supported by an INDAM grant.
}

\begin{abstract}
We classify Fano fivefolds of index two which are blow-ups
of smooth manifolds along a smooth center.
\end{abstract}

\section{Introduction}
A smooth complex projective variety $X$ is called {\it Fano} if its
anticanonical bundle $-K_X$ is ample; the {\it index} $r_X$ of $X$ 
is the largest natural number $m$ such that $-K_X=mH$ for some (ample) 
divisor $H$ on $X$, while the {\it pseudoindex} $i_X$ is the minimum 
anticanonical degree of rational curves on $X$.\par
By a theorem of Kobayashi and Ochiai \cite{KOc}, $r_X = \dim X +1$ if and only if 
$(X,L) \simeq (\proj^{\,\dim X}, \Ol_{\proj}(1))$, and $r_X= \dim X$ if and only if 
$(X,L) \simeq (\quadr^{\,\dim X}, \Ol_{\quadr}(1))$, where $\quadr^{\,\dim X}$ is a quadric 
hypersurface in $\proj^{\,\dim X+1}$.
Fano manifolds of index equal to $\dim X-1$ and to $\dim X-2$, which are called {\it del Pezzo} 
and {\it Mukai} manifolds respectively, have been classified, mainly by Fujita, Mukai and Mella
(see \cite{Fuj,Mu,Me}).
In case of index equal to $\dim X -3$, the classification has been completed for Fano manifolds of 
Picard number $\rho_X$ greater than one and dimension greater or equal than six 
(see \cite{Wifala}).\par
For Fano manifolds of dimension five and index two it was proved in \cite{ACO}
that the Picard number is less than or equal to five, equality holding only for
a product of five copies of $\pu$. 
Then, in \cite{CO}, the structure of the possible Mori cones of curves of those manifolds,
i.e., the number and type of their extremal contractions, was described.
A first step in going from the table of the cones
given in \cite{CO} to the actual classification of Fano fivefolds of
index two has been done in \cite{NO}, where ruled Fano fivefolds of index two, i.e., 
fivefolds of index two
with a $\pu$-bundle structure over a smooth fourfold, were classified.\par
\ms
In this paper we classify Fano fivefolds of index two which are blow-ups of smooth
manifolds along smooth centers. In section \ref{start} we recall the structure of the 
cones of curves of these manifolds, as described in \cite{CO}, 
and we summarize the known results. Using previous results we are reduced to the following cases:
\begin{enumerate} 
\item[] $\rho_X=2$ and the two extremal rays of $\cone(X)$
correspond respectively to the blow-up of a smooth variety $X'$ along a smooth surface $S$ and to a
fiber type contraction $\vartheta:X \to Y$.
\item[] $\rho_X=3$. In this case $\cone(X)$ has three extremal rays: one of them is associated to 
the blow-up of a smooth variety along a smooth surface, one corresponds to a fiber type 
contraction, and the last one is associated
either to another blow-up contraction or to another fiber type contraction.
\end{enumerate}

\ms
The hardest case, which is the heart of the paper and is dealt with in Section \ref{rho2}, is 
when $\rho_X=2$.
In this case it is easy to show that the pseudoindex of $X'$ is equal either to six or to four: 
if $i_{X'}=6$ then
$X' \simeq \proj^5$ by results in \cite{Kepn}, and the classification of $S$
follows observing that $S$ cannot have proper trisecants.
In case $i_{X'}=4$ we prove that also $r_{X'}=4$, i.e., that $X'$ is a del Pezzo manifold and 
that $S$ is a del Pezzo surface.
The classification of $(X',S)$ then follows studying the possible conormal bundles $N^*_{S/X'}$.

\ms
In Section \ref{rho3} we study the case $\rho_X=3$; apart from one case,
the target of the birational contraction is a Fano manifold, which is either a product
with $\pu$ as a factor or a $\proj^3$-bundle over a surface; the classification of the
center follows.

\ms
Our results are summarized in the following

\begin{theorem} \label{main}
Let $X$ be a Fano fivefold of index two which is the blow-up
of a smooth variety $X'$ along a smooth subvariety $S.$ 
Then $(X',S)$ is as in the following table, where,
in the last column, $F$ denotes a fiber type extremal ray$,$ $D_i$
denotes a birational extremal ray whose associated contraction contracts a divisor
to an $i$-dimensional variety and $S$ denotes a ray whose associated contraction is small.
\end{theorem}

{\renewcommand{\baselinestretch}{1.5} \scriptsize
\begin{center}
\begin{tabular}{|c||l|l|l|c|}\hline
$\rho_X$  & No. & $X'$ & $S$ & $NE(X)$\\
\hline\hline
2 & \rm(a1) & $\proj^5$ & a point & $\langle F,D_0 \rangle$ \\
\hline\hline
& \rm(b1)  & $\proj^5$ & a linear $\pd$  & $\langle F,D_2 \rangle$ \\
\hline 
& \rm(b2)  & $\proj^5$ & the complete intersection of three quadrics & $\langle F,D_2 \rangle$ \\
\hline 
& \rm(b3)  & $\proj^5$ & $\pu \times \pu$ embedded by $\Ol(1,2)$ & $\langle F,D_2 \rangle$ \\
\hline 
& \rm(b4)  & $\proj^5$ &  ${\mathbb F}_2$ embedded by $C_0+3f$ & $\langle F,D_2 \rangle$  \\
\hline 
& \rm(b5)  & $\proj^5$ & \begin{tabular}{l} the blow-up of $\pd$ in four points $x_1, \ldots, x_4$ such that \\
the line bundle $\Ol_\pd(3) - \sum E_i$ is very ample \end{tabular} & $\langle F,D_2 \rangle$ \\
\hline 
& \rm(b6)  & $\proj^5$ & \begin{tabular}{l} the blow-up of $\pd$ in seven points $x_0, \ldots, x_6$ such that \\
the line bundle $\Ol_\pd(4) -2E_0 -\sum_{i=1}^6 E_i$ is very ample \end{tabular} 
& $\langle F,D_2 \rangle$ \\
\hline 
& \rm(b7)  & $V_d$ (*)
& the complete intersection of three general members of $|\Ol_{V_d}(1)|$ &
 $\langle F,D_2 \rangle$ \\
\hline 
& \rm(b8)  & $V_3$ & $\pd$ with $(\Ol_{V_3}(1))_{|\pd} \simeq \Ol_\pd(1)$  & $\langle F,D_2 \rangle$ \\
\hline 
& \rm(b9)  & $V_4$ & $\pd$ with $(\Ol_{V_4}(1))_{|\pd} \simeq \Ol_\pd(1)$ & $\langle F,D_2 \rangle$ \\
\hline 
& \rm(b10)  & $V_4$  & $\quadr^2$ with $(\Ol_{V_4}(1))_{|\quadr} \simeq \Ol_\quadr(1)$   &
$\langle F,D_2 \rangle$ \\
\hline 
& \rm(b11)  & $V_5$ & a plane of bidegree  $(1,0)$ (**)  & $\langle F,D_2 \rangle$ \\
\hline 
& \rm(b12)  & $V_5$ & a quadric of bidegree $(1,1)$  & $\langle F,D_2 \rangle$ \\
\hline 
& \rm(b13)  & $V_5$ & a surface $\mathbb F_1$ of bidegree $(2,1)$ not contained in a $\G(1,3)$  & 
$\langle F,D_2 \rangle$ \\
\hline 
\hline
& \rm(c1)  & $\proj^5$ & a Veronese surface  & $\langle D_2,D_2 \rangle$ \\
\hline 
& \rm(c2)  & $\proj^5$ & $\mathbb F_1$ embedded by $C_0 + 2f$  & $\langle D_2,D_2 \rangle$ \\
\hline 
& \rm(c3)  & $V_5$ & a plane of bidegree $(0,1)$  & $\langle D_2,D_2 \rangle$ \\
\hline 
\hline 
& \rm(d1)  & $\proj^5$ & $\quadr^2$ with $(\Ol_{\proj}(1))_{|\quadr} \simeq \Ol_{\quadr(1)}$ & 
$\langle D_2,S \rangle$ \\
\hline \hline 
3 & \rm(e1)  & $\pu \times \quadr^4$ & $\pu \times l$ with $l$ a line in $\quadr^4$  & $\langle F,F,D_2 \rangle$ \\
\hline 
& \rm(e2)  &  $\pu \times \quadr^4$ & $\pu \times \Gamma$ with $\Gamma \subset \quadr^4$ a conic not 
contained in a plane $\Pi  \subset \quadr^4$ & $\langle F,F,D_2 \rangle$ \\
\hline
& \rm(e3)  & $X' \in |\Ol_{\pd \times \pf}(1,1)|$ &  $\pd$, a fiber of the projection
$X' \to \pf$ & $\langle F,F,D_2 \rangle$ \\ 
\hline 
& \rm(e4)  & $X' \in |\Ol_{\pd \times \pf}(1,1)|$ & \begin{tabular}{l} $\mathbb F_1$, the complete intersection
of $X'$ and three general \\  members of the linear system $|\Ol_{\pd \times \pf}(0,1)|$ 
\end{tabular} & $\langle F,F,D_2 \rangle$ \\ 
\hline \hline 
& \rm(f1)  & $\proj_\pd(\Ol \oplus \Ol(1)^{\oplus 3})$ &   $\pd$, a section corresponding
to the surjection $\Ol\oplus \Ol(1)^{\oplus 3} \to \Ol$ & $\langle F,D_2,D_2 \rangle$ \\ 
\hline
& \rm(f2)  & $\Bl_\pi(\proj^5)$ (***)& $\pd$, a non trivial fiber
of $\Bl_{\pi}(\proj^5) \to \proj^5$  & 
$\langle F,D_2,D_2 \rangle$ \\ 
\hline
& \rm(f3)  & $\Bl_p(\proj^5)$ & $\mathbb F_1$, the strict transform of a plane in $\proj^5$ through $p$  & 
$\langle F,D_2,D_2 \rangle$ \\  
\hline
& \rm(f4)  & $\Bl_\pi(\proj^5)$ & $\pd$, the strict transform of a plane in $\proj^5$ not meeting $\pi$  & 
$\langle F,D_2,D_2 \rangle$ \\ 
\hline \hline 
4 & \rm(g1)  & $\pu \times \pu \times \pt$ & $\pu \times \pu \times \{p\}$  & $\langle F,F,F,D_2 \rangle$ \\
\hline 
\end{tabular}
\end{center}}
\renewcommand{\baselinestretch}{1.2} 
\footnotetext{(*) $V_d$ denotes a del Pezzo fivefold of degree $d$.} 
\footnotetext{(**)  $V_5$ is a hyperplane section of $\G(1,4)$. The bidegree of $S$
is the bidegree of $S$ in $\G(1,4)$.} 
\footnotetext{(***) $\Bl_\pi(\proj^5)$ (resp. $\Bl_p(\proj^5)$) denotes the blow-up of
$\proj^5$ along a $2$-plane $\pi$ (resp. along a point $p$).} 
\ms

\noindent
In \cite{AOlong}, Fano manifolds $X$ obtained by blowing up a smooth variety $Y$ along a center
$T$ of dimension $\dim T \le i_X-1$ were classified; the results in this paper show
that the case $\dim T = i_X$ will be far more complicated.


\section{Preliminaries}

\subsection{Fano-Mori contractions and rational curves}

Let $X$ be a smooth Fano variety of dimension $n$ and $K_X$  its canonical divisor.
By Mori's {\it Cone Theorem} the cone $\cone(X)$ of effective 1-cycles, which is contained in
the $\mathbb R$-vector space $N_1(X)$ of 1-cyles modulo numerical equivalence,
is polyhedral; a face $\tau$ of $\cone(X)$ is called an 
{\it extremal face}  
and an extremal face of dimension one is called an {\it extremal ray}.
To every extremal face $\tau$ one can associate a morphism $\f:X \to Z$ with connected fibers 
onto a normal variety; the morphism $\f$ contracts those curves whose numerical class lies in $\tau$, 
and is usually called the {\it Fano-Mori contraction}
(or the {\it extremal contraction}) associated to the face $\tau$.
A Cartier divisor $D$ such that $D = \f ^*A$ for an ample divisor $A$ on $Z$
is called a {\it supporting divisor} of the map $\f$ (or of the face $\tau$).\\
An extremal ray $R$ is called {\it numerically effective}, or of
{\it fiber type}, if $\dim Z < \dim X$, otherwise the ray is {\it non nef} or {\it birational}.
We usually denote with $E = E(\f):= \{ x \in X\ |\  \dim\f^{-1} (\f(x)) > 0\}$
the {\it exceptional locus} of $\f$; if $\f$ is of fiber type then of course $E=X$.
If the exceptional locus of a birational ray $R$ has codimension one,
the ray and the associated contraction are called {\it divisorial}, otherwise they are called {\it small}.

\begin{definition} \label{defcontrazioni}
An elementary fiber type extremal contraction 
$\f: X \to Z$ is called a {\it scroll} (resp. a {\it quadric fibration}) 
if there exists a $\f$-ample line bundle $L \in \pic(X)$ such that
$K_X+(\dim X-\dim Z+1)L$ (resp. $K_X+(\dim X-\dim Z)L$)
is a supporting divisor of $\f$.
An elementary fiber type extremal contraction $\f:X \to Z$ onto a smooth
variety $Y$ is called a $\proj${\it-bundle} if there exists a vector bundle $\E$ of rank 
$\dim X-\dim Z+1$ on $Z$ such that 
$X \simeq \proj_Z(\E)$;
every equidimensional scroll is a $\proj$-bundle by \cite[Lemma 2.12]{Fu1}.
\end{definition}

\begin{definition} Let $\ratcurves$ be the normalized space of rational curves in $X$ in the sense of
\cite{Kob};  a {\it family of rational curves} will be an irreducible component
$V \subset \ratcurves$.
Given a rational curve $f:\pu \to X$ we  call a {\it family of
deformations} of \linebreak $f$ any irreducible component $V \subset
\ratcurves$ containing the equivalence class of $f$.
\end{definition}

We define $\loc(V)$ to be the subset of $X$ of points through which there
is a curve parametrized by $V$; we say that $V$ is a {\it dominating family}
if $\overline{\loc(V)}=X$.
Moreover, for every point $x \in \loc(V)$, we will denote by $V_x$ the 
subscheme of $V$ parametrizing rational curves passing through $x$.\par

\begin{definition} 
Let $V$ be a family of rational curves on $X$. We say that $V$ is {\it unsplit} if it is proper and that
$V$ is {\it locally unsplit} if every component of $V_x$ is proper for the general $x \in \loc(V)$.
\end{definition}

\begin{proposition}\cite[IV.2.6]{Kob}\label{iowifam} 
Let $X$ be a smooth projective variety, $V$ a family of rational curves
and $x \in \loc(V)$ such that every component of $V_x$ is proper. Then

      \item[(a)] $\dim X  -K_X \cdot V \le \dim \loc(V)+\dim \loc(V_x) +1$;
      \item[(b)] $-K_X \cdot V \le \dim \loc(V_x)+1$.

\end{proposition}

In case $V$ is the unsplit family of deformations of a minimal extremal
rational curve, Proposition \ref{iowifam} gives the {\it fiber locus inequality}:

\begin{proposition}{\rm\cite{Io1,Wicon}}\label{fiberlocus} 
Let $\f$ be a Fano-Mori contraction
of $X$ and  $E$ its exceptional locus$.$
Let $F$ be an irreducible component of a {\rm (non trivial)} fiber of $\f$. Then
$$\dim E + \dim F \geq \dim X + l -1$$
where $l =  \min \{ -K_X \cdot C\ |\  C \textrm{~is a rational curve in~} F\}.$\\
If $\f$ is the contraction of a ray $R,$ then $l$ is called the {\it length of the ray}.
\end{proposition}

\begin{definition}\label{CF}
We define a {\it Chow family of rational curves} $\V$ to be an irreducible 
component of  $\textrm{Chow}(X)$ parametrizing rational and connected 1-cycles.\\
If $V$ is a family of rational curves$,$ the closure of the image of
$V$ in $\textrm{Chow}(X)$ is called the {\it Chow family associated to} $V$.
\end{definition}

\begin{definition}
Let $X$ be a smooth variety$,$ $\V^1, \dots, \V^k$ Chow families of rational curves
on $X$ and $Y$ a subset of $X.$ We denote by $\loc(\V^1, \dots, \V^k)_Y$ the set of points $x \in X$ 
that can be joined to $Y$ by a connected chain of $k$ cycles belonging \underline{respectively} to the families 
$\V^1, \dots, \V^k.$ We denote by $\cloc_m(\V^1, \dots, \V^k)_Y$ the set of points $x \in X$ that can be joined 
to $Y$ by a connected chain  of at most $m$ cycles belonging to the families $\V^1, \dots, \V^k$.
\end{definition}

\begin{definition} 
Let $V^1, \dots, V^k$ be unsplit families on $X$.
We will say that $V^1, \dots, V^k$ are {\it numerically independent} if their numerical classes
$[V^1], \dots ,[V^k]$ are linearly independent in the vector space $N_1(X)$.
When moreover $C \subset X$ is a curve, we will say that $V^1, \dots, V^k$ are numerically 
independent from $C$ if the class of $C$ in $N_1(X)$ is not contained in the
vector subspace generated by $[V^1], \dots ,[V^k]$.   
\end{definition}

\begin{lemma} \label{locy}{\rm\cite[Lemma 5.4]{ACO}}
Let $Y \subset X$ be a closed subset and $V$ an unsplit family.
Assume that curves contained in $Y$ are numerically independent from curves in $V,$ and that
$Y \cap \loc(V) \not= \emptyset$.
Then for a general $y \in Y \cap \loc(V)$

      \item[\rm(a)] $\dim \loc(V)_Y \ge \dim (Y \cap \loc(V)) + \dim \loc(V_y);$
      \item[\rm(b)] $\dim \loc(V)_Y \ge \dim Y -K_X \cdot V - 1$.

\item Moreover, if $V^1, \dots, V^k$ are numerically independent
unsplit families such that curves contained in $Y$ are numerically independent
from curves in $V^1, \dots, V^k$, then either \linebreak
$\loc(V^1, \ldots, V^k)_Y=\emptyset$ or

      \item[\rm(c)] $\dim \loc(V^1, \ldots, V^k)_Y \ge \dim Y +\sum (-K_X \cdot V^i) -k$.

\end{lemma}

\begin{definition}
We define on $X$ a relation of {\it rational connectedness with respect to $\V^1, \dots, \V^k$}
 in the following way$:$ $x$ and $y$ are in rc$(\V^1,\dots,\V^k)$-relation if there
exists a chain of rational curves in $\V^1, \dots ,\V^k$ which joins $x$ and $y,$ i.e.
if $y \in \cloc_m(\V^1, \dots, \V^k)_x$ for some $m$. If all the points of $X$
are in rc$(\V^1,\dots,\V^k)$-relation we say that $X$ is {\it rc$(\V^1,\dots,\V^k)$-connected}.
\end{definition}

To the rc$(\V^1,\dots,\V^k)$-relation we can associate a fibration, at least on an open subset.
(see \cite[IV.4.16]{Kob}); we will call it rc$(\V^1,\dots,\V^k)$-fibration.

\begin{definition}
Let $\V$ be the Chow family associated to a family of rational curves $V$. We say that
$V$ is {\it quasi-unsplit} if every component of any reducible cycle in $\V$ is
numerically proportional to $V$.
\end{definition}

{\it Notation}: Let $T$ be a subset of $X$. 
We write $N_1^X(T)=\morespan{V^1}{V^k}$ if the numerical class in $X$ of every curve 
$C \subset T$ can be written as $[C]= \sum_i a_i [C_i]$, with $a_i \in \quadr$ and 
$C_i \in V^i$. We write $\cone^X(T)=\morespan{V^1}{V^k}$  (or $\cone^X(T)=\morespan{R_1}{R_k}$) 
if the numerical class in $X$ of every curve $C \subset T$ can
be written as $[C]= \sum_i a_i [C_i]$, with $a_i \in \quadr_{\ge 0}$ 
and $C_i \in V^i$ (or $[C_i]$ in $R_i$).\par

\begin{proposition}{\rm \cite[Corollary 4.2]{ACO}, \cite[Corollary 2.23]{CO}} \label{coneS}
Let $V$ be a family of rational curves and $x$ a point in $\loc(V)$.

\item[\rm(a)]  If $V$ is quasi-unsplit, then 
$\cone^X(\cloc_m(V)_x)=\onespan{V}$ for every $m \ge 1;$
\item[\rm(b)] 
if $V_x$ is unsplit, then $\cone^X(\loc(V_x))=\onespan{V}$.

\item Moreover, if $\tau$ is an extremal face of $\cone(X),$ $F$ is a fiber
of the associated contraction and $V$ is unsplit and independent from $\tau,$
then

\item[\rm(c)] $\cone^X(\cloc_m(V)_F) = \langle \tau, [V] \rangle$ for every $m \ge 1$.

\end{proposition}



\subsection{Fano bundles}\label{fb}

\begin{definition}
Let $\E$ be a vector bundle on a smooth complex projective variety $Z$. We say that $\E$ is a {\it Fano bundle} if 
$X=\proj_Z(\E)$ is a Fano manifold. By \cite[Theorem 1.6]{SzW2}\label{fanooverfano}
if $\E$ is a Fano bundle over $Z$ then $Z$ is a Fano manifold.
\end{definition}

M. Szurek and J. Wi\'sniewski have classified Fano bundles over $\pd$ (\cite{SzW1, SzW3}) and 
Fano bundles of rank two on surfaces \cite{SzW3}.
What follows is a characterization of Fano bundles of rank $r \ge 2$ over del Pezzo surfaces, 
which generalizes some results in \cite{SzW3}.

\begin{proposition} \label{fbodps}
Let $S_k$ be a del Pezzo surface obtained by blowing up $k >0$ points in $\pd,$ and let
$\E$ be a Fano bundle of rank $r \ge 2$ over $S_k;$ then$,$ up to twist $\E$ with a suitable line bundle$,$
the pair $(S_k,\E)$ is one of the following$:$
\item[\rm(i)] $(S_k,\oplus \Ol^{\oplus r});$
\item[\rm(ii)] $(S_1, \theta^*(\Ol_\pd(1) \oplus \Ol_\pd^{\oplus (r-1)}));$
\item[\rm(iii)] $(S_1, \theta^*(T\pd(-1) \oplus \Ol_\pd^{\oplus (r-2)})),$
\item where $\theta:S_1 \to \pd$ is the blow-up of $\pd$ at one point.
\end{proposition}

\begin{proof} 
Let $\E$ be a Fano bundle of rank $r \ge 2$ over $S_k$ and let $X=\proj_{S_k}(\E)$;
by \cite[Proposition 3.4]{NO} there is a one-to-one correspondence between
the extremal rays of $\cone(S_k)$ and the extremal rays of $\cone(X)$ spanning a 
two-dimensional face with the ray $R_\E$ corresponding to the projection
$p:X \to S_k$.\\
Let $R_{\theta_1} \subset \cone(S_k)$ be an extremal ray of $S_k$ associated to a blow-up
${\theta_1}: S_k \to S_{k-1}$, and call $E_{\theta_1}$ the exceptional divisor of ${\theta_1}$; 
let $R_{\vartheta_1}$ be the corresponding ray in
$\cone(X)$, with associated extremal contraction ${\vartheta_1}:X \to X_1$.\\
By \cite[Lemma 3.5]{NO} ${\vartheta_1}$ is birational and has one-dimensional fibers,
hence by \cite[Theorem 5.2]{AOsperays} we have that $X_1$ is smooth and ${\vartheta_1}$ 
is the blow-up of a smooth subvariety of codimension two in $X_1$; moreover, by
\cite[Lemma 3.5]{NO} and reasons of dimension, $\Exc(R_{\vartheta_1})= p^{-1}(E_{\theta_1})$.\\
The divisor $E_{\vartheta_1}:=\Exc(R_{\vartheta_1}) $ has two projective bundle structures: 
a $\proj^1$ bundle structure over 
the center of the blow-up and a $\proj^{r-1}$-bundle structure over $E_{\theta_1}$;
by \cite[Main theorem]{Sa}  we have that $E_{\vartheta_1} \simeq \pu \times \proj^{r-1}$.\\
It follows that $\E_{|E_{\theta_1}} \simeq \Ol^{\oplus r}$, hence by \cite[Lemma 2.9]{AOext} 
there exists a vector bundle of rank $r$ on $S_{k-1}$ such that $\E=\theta_1^*\E_1$. 
It is now easy to prove that the induced map $\proj_{S_k}(\theta_1^*\E_1)= X \longrightarrow
\proj_{S_{k-1}}(\E_1)$ is nothing but ${\vartheta_1}$, hence $X_1= \proj_{S_{k-1}}(\E_1)$.\\
Since $\cone(E_{\vartheta_1})=\langle R_\E, R_{\vartheta_1}\rangle$, the divisor $E_{\vartheta_1}$ 
cannot contain the exceptional locus of another extremal ray of $X$; it follows that
$X_1$ is a Fano manifold  by \cite[Proposition 3.4]{Wicon}.\par

\ms
We iterate the argument $k$ times, until we find a Fano bundle $\E_k$ over $\pd$
such that, denoted by $\theta$ and $\vartheta$ the composition of the contractions $\theta_i$ and 
$\vartheta_i$ respectively, $\E=\theta^*\E_k$. We have a commutative diagram 
$$\xymatrix@=35pt{
\proj_{S_k}(\E)=X \ar[r]^{\vartheta} \ar[d]_{p}  & X_k=\proj_\pd(\E_k) \ar[d]^{p_k} \\
S_k \ar[r]_{\theta} & \pd }$$

Up to considering the tensor product of $\E_k$ with a suitable line bundle, we can assume that
$0 \le c_1(\E_k) \le r-1$; by \cite[Proposition 2.2]{SzW1} we have that $\E_k$ is nef.\par

\ms
Let $l$ be a line in $\pd$; the restriction of $\E_k$ to $l$ decomposes as a sum of
nonnegative line bundles, hence we can write $(\E_{k})_{|l} \simeq \oplus_{i=0}^{r-1} \Ol(a_i)$,
with $a_0=0$ and $a_i \ge 0$.
Let $\tl l$ be the strict transform of $l$ in $S_k$; since $\theta_{|\tl l}: \tl l \to l$
is an isomorphism we have $\E_{|\tl l} \simeq (\E_{k})_{|l}$; let $C_0 \subset X$
be a section of $p$ over $\tl l$ corresponding to a surjection $\E_{|\tl l} \to \Ol \to 0$; 
we have
\begin{equation}\label{conto}
0<-K_X \cdot C_0 = ra_0 -K_{S_k} \cdot \tl l -\sum_{i=0}^{r-1} a_i = -K_{S_k} \cdot \tl l - c_1(\E_k).
\end{equation}
Now if $l$ passes through a point blown up by $\theta$, by equation (\ref{conto}) we have 
$c_1(\E_k) \le 1$.
In this case, by the classification in \cite{SzW1}, either $\E_k$ is trivial, or
$\E_k \simeq \Ol(1) \oplus \Ol^{\oplus(r-1)}$, or $\E_k \simeq T\pd(-1) \oplus \Ol^{\oplus(r-2)}$.

Assume that $k \ge 2$ and let $l$ be a line in $\pd$ joining two 
of the blown-up points; again by equation (\ref{conto}) we have $c_1(\E_k)=0$,
so only the first case occurs.
\end{proof}

\begin{proposition}\label{fanoquadr}
Let $\E$ be a Fano bundle of rank $r \ge 2$ over $\pu \times \pu;$ then$,$ up to twist
$\E$ with a suitable line bundle$,$  $\E$ is one of the following$:$

\item[\rm(i)] $\Ol^{\oplus r};$ 
\item[\rm(ii)] $\Ol(1,0) \oplus \Ol^{\oplus(r-1)};$ 
\item[\rm(iii)] $\Ol(1,1) \oplus \Ol^{\oplus(r-1)};$
\item[\rm(iv)] $\Ol^{\oplus (r-2)} \oplus \Ol(1,0) \oplus \Ol(0,1);$
\item[\rm(v)] a vector bundle fitting in the exact sequence 
$0 \to \Ol(-1,-1) \to \Ol^{\oplus (r+1)} \to \E \to 0$.
\item
In all cases the cone of curves of $X=\proj(\E)$ is generated by the ray corresponding to the bundle projection
and by other two extremal rays; in case (i) the other rays are of fiber type, in case (ii) one of them 
is of fiber type and the other correspond to a smooth blow-up while in cases (iii)-(v) both the other rays
correspond to smooth blow-ups.
\end{proposition}

\begin{proof} We will show the result by induction on $r$, the case $r=2$ having been 
established in  \cite[Main Theorem]{SzW3}.
Let $X=\proj(\E)$; first of all we prove that $\cone(X)$ is generated by three extremal rays.
Let $R_\E \subset \cone(X)$ be the extremal ray corresponding to the projection $p:X \to \pu \times \pu$;
since $\rho_X=3$ it is enough to prove that any other extremal ray of $\cone(X)$ lies in a two-dimensional
face with $R_\E$.\par
Let $R_{\vartheta}$ be another extremal ray of $X$ with associated contraction $\vartheta$ and let $F$ be a nontrivial
fiber of $\vartheta$. We claim that $\dim F=1$: in fact, 
since curves contained in $F$ are not contracted by $p$, we have $\dim F \le 2$, and,
if $\dim F = 2$, we would have $X = p^{-1}(p(F))$ and $\cone(X)=\langle R,R_\E \rangle$
by Proposition \ref{coneS} (c), against the fact that $\rho_X=3$.
In particular, by Proposition \ref{fiberlocus}, $\vartheta$ cannot be a small contraction.\par
\ms
Let $V_\vartheta$ be a family of rational curves of minimal degree (with respect to some fixed ample line bundle) 
among the families which dominate the exceptional locus of $R_{\vartheta}$ and whose class 
is in $R_{\vartheta}$. Such a family is quasi-unsplit by the extremality of $R_{\vartheta}$
and locally unsplit by the assumptions on its degree.\\
We claim that $V_\vartheta$ is horizontal and dominating with respect to $p$.
This is clear if the contraction $\vartheta$ associated to $R_{\vartheta}$ is of fiber type.
Assume that $\vartheta$ is divisorial, with exceptional locus $E$: we cannot have
$E \cdot R_\E =0$, otherwise $E=p^*D$ for some effective divisor $D$ in $\pu \times \pu$;
but every effective divisor on $\pu \times \pu$ is nef and so $E$ would be nef,
against the fact that $E \cdot R_{\vartheta} <0$.
It follows that $E \cdot R_\E >0$, so $E$ dominates $\pu \times \pu$ and thus $V_\vartheta$
is horizontal and dominating with respect to $p$, and the claim is proved.\\
We can now apply \cite[Lemma 2.4]{CO} and conclude that $[V_\vartheta]$ and $R_\E$ lie
in a two-dimensional extremal face of $\cone(X)$.\par
\ms
We have thus proved that every extremal ray different from $R_\E$ lies in a two-dimensional
face with $R_\E$; therefore $\cone(X)$ is generated by three extremal rays.
We will call $R_{\vartheta_1}$ and $R_{\vartheta_2}$ the two rays different from $R_\E$, i.e.,
$\cone(X)=\langle R_\E, R_{\vartheta_1}, R_{\vartheta_2}\rangle$.\par
\ms
By \cite[Proposition 3.4]{NO}, for every $i=1,2$ we have a commutative diagram 
$$\xymatrix@=35pt{
X \ar[r]^{\vartheta_i} \ar[dr]^{\psi_i} \ar[d]_p  & Z_i \ar[d] \\
\pu \times \pu \ar[r]_{\theta_i} & \pu }$$
where $\psi_i$ is the contraction of the face of $\cone(X)$ 
spanned by $R_\E$ and $R_{\vartheta_i}$. 

Let $x \in \pu$ and let $f_x^i$ be the fiber of $\theta_i$ over $x$;
since we can factor $\psi_i$  as $\psi_i= \theta_i \circ p$, the fiber 
of $\psi_i$ over $x$ is $\proj(\E_{|f^i_x})$. By the smoothness of $\psi_i$
and adjunction, $\proj(\E_{|f^i_x})$ is a Fano manifold, hence 
either $\E_{|f^i_x} \simeq  \Ol(a)^{\oplus r}$ or $\E_{|f^i_x} \simeq \Ol(a+1) \oplus \Ol(a)^{\oplus (r-1)}$.\\
Since the degree of $\E$ does not change as $x$ varies in $\pu$ we have
that, for a fixed $i=1,2$, the splitting type of $\E$ along the fibers of $\theta_i$
is constantly $(a,\dots,a)$ or $(a+1,a, \dots,a)$.\
Up to twist $\E$ with a line bundle we can assume that
its  splitting type along the fibers of $\theta_i$ 
is constantly $(0,\dots,0)$ or $(1,0, \dots,0)$
\par
\ms
If for some $i=1,2$ the splitting type of $\E$ on the fibers of $\theta_i$ is $(0, \dots,0)$ 
then $\E \simeq \theta_i^*\E'$, with $\E'$ a vector bundle on $\pu$; hence $\E$ 
is decomposable and we are in case (i) or (ii).\par
\ms
Assume now that the splitting type of $\E$ on the fibers of $\theta_i$ 
is $(1,0,\dots,0)$ for $i=1$ and $i=2$, and thus $c_1(\E)=(1,1)$.
We claim that in this case the contractions $\vartheta_i:X \to Z_i$ are birational.\\
Assume by contradiction that for some $i$, say $i=1$, the contraction $\vartheta_1$ is of fiber type.
Let $x \in \pu$ be a general point; the fiber of $Z_1 \to \pu$ has dimension strictly smaller
than the dimension of $\psi_1^{-1}(x)$. It follows that both the restrictions of $\vartheta_1$
and $p$ to $\psi_1^{-1}(x)$ are of fiber type, yet $\psi_1^{-1}(x)\simeq \Bl_{\proj^{r-2}}(\proj^{r})$,
so it has only one fiber type contraction.\par

\ms
We have already proved that the nontrivial fibers of the contractions $\vartheta_i$ are one dimensional,
hence for every $i=1,2$ the variety $Z_i$ is smooth and $\vartheta_i$ is the blow-up of a smooth 
subvariety of codimension two in $Z_i$ by \cite[Theorem 5.2]{AOsperays}.\\
Consider one of the birational contractions of $X$, say $\vartheta_1: X \longrightarrow Z_1$, and
let $E_1$ be its exceptional locus. 
For every fiber $f_x$ of $\theta_1$ the restriction of $E_1$ to $\proj_{f_x}(\E_{|f_x})$ is
a non nef divisor, hence it is the exceptional divisor of the contraction
$\proj_{f_x}(\E_{|f_x}) \to \proj^r$. In particular $E_1 \cdot R_\E=1$ and  $E_1$ does not contain any fiber of $p$.
By \cite[Lemma 2.12]{Fu1} the restriction of $p$ makes $E_1$ a projective bundle over $\pu \times \pu$, that is
$E_1=\proj_{\pu \times \pu}(\E')$ with $\E'$ a rank $r-1$ vector bundle over $\pu \times \pu$.
We will now split the proof in two cases, depending on the sign of the intersection number of $E_1$
with $R_{\vartheta_2}$.\par
\ms
{\bf Case 1.} \quad $E_1 \cdot R_{\vartheta_2} \le 0$.\par
\ms
In this case the line bundle $-K_X -E_1$ is ample on $X$; therefore its restriction
to $E_1$ is ample, $E_1$ is a Fano manifold and $\E'$ is
a Fano bundle of rank $r-1$ over $\pu \times \pu$.\\
Note also that $E_1$ has a fiber type contraction different from the bundle projection onto
$\pu \times \pu$, coming by the blow-up contraction $\vartheta_1$, 
so, by induction, either $\E'$ is trivial or $\E' \simeq \Ol(1,0) \oplus \Ol^{\oplus (r-2)}$. 
The injection 
$\proj_{\pu \times \pu}(\E') \hookrightarrow \proj_{\pu \times \pu}(\E)$
gives an exact sequence of bundles on $\pu \times \pu$
$$\xymatrix{0 \ar[r] & \Ol(a,b) \ar[r] & \E \ar[r] & \E' \ar[r]& 0},$$
with $E_1= \xi_\E +p^*\Ol(-a,-b)$. Computing the intersection numbers of $E_1$ with $R_{\vartheta_1}$
and $R_{\vartheta_2}$ and recalling the splitting type of $\E$ we have the following
possibilities:
$$\xymatrix{0 \ar[r] & \Ol(0,1) \ar[r] & \E \ar[r] &  \Ol^{\oplus(r-2)} \oplus \Ol(1,0) \ar[r]& 0};$$
$$\xymatrix{0 \ar[r] & \Ol(1,1) \ar[r] & \E \ar[r] &  \Ol^{\oplus(r-1)}  \ar[r]& 0}.$$

Both these sequences split, so we are in cases (iii) or (iv).\par
\ms
{\bf Case 2.} \quad $E_1 \cdot R_{\vartheta_2} > 0$.\par
\ms
By \cite[Proposition 3.4]{Wicon} $Z_1$ is a Fano manifold. $Z_1$ has a fiber type elementary contraction onto $\pu$.
For a general $x \in \pu$ the fiber $\psi_1^{-1}(x)=\proj(\E_{|f^i_x})$ is isomorphic to $\Bl_{\proj^{r-2}}(\proj^{r})$,
hence the fiber of $Z_1 \to \pu$ over $x$ is isomorphic to $\proj^r$. It follows that $Z_1$
has a projective bundle structure over $\pu$ (cfr. \cite[Lemma 2.17]{NO}),
so either $Z_1 \simeq \pu \times \proj^r$ or $Z_1 \simeq \Bl_{\proj^{r-1}}(\proj^{r+1})$.\\
The second case cannot happen: in fact, let $\psi:X \to \proj^{r+1}$
be the contraction of the face spanned by $R_{\vartheta_1}$ and $R_{\vartheta_2}$. 
Denoting
by $E$ the exceptional divisor of the contraction $Z_1 \to \proj^r$, by $\tl E$ its
strict transform in $X$, and applying twice
the canonical bundle formula for blow-ups  we have
$$K_X= \vartheta_1^*K_{Z_1}+E_1= \psi^*K_{\proj^{r+1}}+\vartheta_1^*E +
 E_1=\psi^*K_{\proj^{r+1}}+\tl E + kE_1.$$
Since $K_X \cdot R_{\vartheta_2} =-1 $ and $\psi^*K_{\proj^{r+1}}\cdot R_{\vartheta_2} =0$
we have $\tl E \cdot R_{\vartheta_2}<0$.
This implies that $\tl E =E_2$, and thus $\tl E \cdot R_{\vartheta_2}=-1$, 
yielding $E_1 \cdot R_{\vartheta_2}=0$, a contradiction.\par
\ms 
Note that the minimal extremal curves contracted by $\vartheta_i$ are the minimal
sections (those corresponding to the trivial summands) of $p:\proj(\E_{|f^i_x}) \to \pu$ along 
the fibers of $\theta_i$; therefore $\xi_\E \cdot R_{\vartheta_i}=0$ for $i=1,2$.
Being trivial on the face spanned by $R_{\vartheta_1}$ and $R_{\vartheta_2}$ and 
positive on $R_\E$ the line bundle $\xi_\E$ is nef.\\
Let $\psi$ be the contraction of the face spanned by $R_{\vartheta_1}$ and $R_{\vartheta_2}$; this contraction
factors through $Z_1 \simeq \pu \times \proj^r$ and therefore is onto $\proj^r$, since it does not
contract curves in $R_\E$.
The line bundle $\xi_\E$ restricts to $\Ol(1)$ on the fibers of $p$, hence $\xi_\E=\psi^*\Ol_{\proj^r}(1)$. 
Therefore $\xi_\E$ (and so $\E$) is spanned and we have an exact sequence on $\pu \times \pu$:
$$\xymatrix{0 \ar[r] & \Ol(a,b) \ar[r] &  \Ol^{\oplus(r+1)} \ar[r] & \E \ar[r]& 0},$$
Computing the first Chern class we have $a=-1,b=-1$ and we are in case (v). In this case
$X=\proj(\E)$ is a divisor in the linear system $\Ol(1,1,1)$ in 
$\proj_{\pu \times \pu}(\Ol^{\oplus(r+1)}) \simeq \pu \times \pu \times \proj^r$.
\end{proof}

\subsection{Surfaces in $\G(1,4)$}\label{sg}

Let $\G(r,n)$ be the Grassmann variety of projective $r$-spaces in $\proj^n$, embedded in $\proj^N$ via the
Pl\"ucker embedding. We will denote a point in $\G(r,n)$ 
by a capital letter, and the corresponding linear space in $\proj^n$ by the same small letter.\par
\medskip
Consider the Schubert cycles $\Omega_1:=\Omega(0,1, \dots, r-1,r+2)$ and $\Omega_2:=\Omega(0,1,\dots,r-2,r,r+1)$;
the cohomology class of a surface $S \subset \G(r,n)$ can be written as $\alpha \Omega_1 +\beta\Omega_2$.
Recalling that the class of an hyperplane section of $\G(r,n)$ is the class of the Schubert cycle
$\Omega_H:=\Omega(n-r-1,n-r,\dots,n-2,n)$ we obtain that the degree of $S$ as a subvariety
of $\proj^N$ is given by 
$$\deg(S)=\alpha \Omega_1\Omega_H^2 +\beta\Omega_2\Omega_H^2 =\alpha +\beta.$$
The integer $\alpha$ is the number of linear spaces parametrized by $S$ which meet a 
general $(n-r-2)$-space in $\proj^n$, as one can see intersecting 
with the Schubert cycle $\Omega(n-r-2,n-r+1,n-r+2, \dots,n)$; it is called the {\it order}  of $S$ 
and denoted by $\ord(S)$.\\
The integer $\beta$ is the number of linear spaces parametrized by $S$ which meet a general $n-r$ space in a line,
as one can see intersecting with the Schubert cycle $\Omega(n-r-1,n-r,n-r+2, \dots,n)$; 
it is called the {\it class} of $S$ and denoted by  $\cl(S)$.

\begin{definition}
The {\it bidegree} of $S$ is the pair $(\ord(S),\cl(S))$. By the discussion above we have that 
$\deg S= \ord(S)+\cl(S)$.
\end{definition}

\begin{remark}\label{rhoplane} 
A $2$-plane $\Lambda_\pi^2$ in $\G(1,4)$ which parametrizes the family of lines which are contained 
in a given $2$-plane $\pi \subset \proj^4$, classically called a {\it $\rho$-plane}, has bidegree $(0,1)$.\\
Moreover, given a point $L \in \G(1,4)$ there exists a line in $\G(1,4)$ joining $\Lambda^2_\pi$ and $L$
if and only if the corresponding line $l \subset \proj^4$ has nonempty intersection with $\pi$.
\end{remark}

\begin{remark}\label{sigmaplane}
The family of lines through a given point $p$ in $\proj^4$ is parametrized by a three-dimensional 
linear space $\Lambda^3_p \subset \G(1,4)$, classically called a $\Sigma$-solid.\\
A two-dimensional linear subspace of a $\Sigma$-solid, classically called a {\it $\sigma$-plane},
parametrizes the family of lines through a given point in $\proj^4$ which lie in a given hyperplane $H$,
and has bidegree $(1,0)$; we will denote it by $\Lambda^2_{p,H}$.\\
Given a $\sigma$-plane $\Lambda^2_{p,H}$ and a point $L \in \G(1,4)$ there exists always a line in 
$\G(1,4)$ joining $\Lambda^2_{p,H}$ and $L$. This is clear if $L$ is contained in the $\Sigma$-solid
$\Lambda^3_p$; otherwise, let $\pi$ be the plane $\subset \proj^4$ 
spanned by $l$ and $p$ and let $q$ be $l \cap H$ if $l \not \in H$ or any point of $l$
if $l \subset H$: the pencil of lines in $\pi$ with center 
$q$ is represented by a line in $\G(1,4)$ passing through $L$ and meeting $\Lambda^2_{p,H}$.
\end{remark}

\begin{example}\label{blowupaplane}
If $\Lambda^2_\pi$ is a $2$-plane of bidegree $(0,1)$ (a $\rho$-plane) then the blow-up
of $\G(1,4)$ along $\Lambda^2_\pi$ is a Fano manifold whose other contraction is the 
blow-up of $\proj^6$ along a cubic threefold contained in a hyperplane (see \cite[Theorem XLI]{SR}).\\
If else  $\Lambda^2_{p,H}$ is a $2$-plane of bidegree $(1,0)$ (a $\sigma$-plane)
the linear system $|\Ol_{\G}(1) \otimes {\cal I}_{\Lambda^2_{p,H}}|$ defines a rational map
${\G} \xymatrix{\ar@{-->}[r] &} \proj^6$  whose image is a quadric cone in $\proj^6$ with zero-dimensional vertex;
the blow-up of $\G(1,4)$ along $\Lambda^2_{p,H}$ is a Fano manifold whose other contraction is of fiber type
onto this quadric cone. This can be checked by direct computation. 
\end{example}

\begin{lemma}\label{zero} 
Let $S$ be a surface in $\G(1,4).$ If \,$\ord(S)=0,$
then $S$ is a plane of bidegree $(0,1),$ while if $\cl(S) = 0,$ 
then $S$ is contained in a $\Sigma$-solid.
\end{lemma}

\begin{proof} 
Let $I \subset \G(1,4) \times \proj^4$ be the incidence variety. Denote
by $p_1:I \to \G(1,4)$ and $p_2:I \to \proj^4$ the projections and let $\loc(S)=p_2(p_1^{-1}(S))$.
If $\ord(S)=0$, then the general line of $\proj^4$ does not meet $\loc(S)$; therefore $\loc(S)$
is two-dimensional. Moreover, since $p_1^{-1}(S)$ is irreducible, also $\loc(S)$ is irreducible.
Therefore $\loc(S)$ is an irreducible surface in $\pf$ which contains a two-parameter
family of lines. It is easy to prove that $\loc(S)$ is a plane,
hence $S$ is the $\rho$-plane which parametrizes the lines of $\loc(S)$.\par
Assume now that $\cl(S)=0$. Since we can identify $\G(1,4)$ with the Grassmannian $\G(2,4)$
of planes in the dual space ${\proj^4}^*$, $S$ can be viewed as a surface
which parametrizes a two-dimensional family of planes in ${\proj^4}^*$.
The duality exchanges order and class, so $S$, as a subvariety of $\G(2,4)$, has order zero, i.e.,
through a general point of ${\proj^4}^*$ there are no planes parametrized by $S$.\\
Denote by $I^* \subset \G(2,4) \times {\proj^4}^*$ the incidence variety, 
by $p_1^*:I^* \to \G(2,4)$ and $p_2^*:I^* \to {\proj^4}^*$ the projections and define
$\loc^*(S)=p_2^*({p_1^*}^{-1}(S))$. Then  $\dim \loc^*(S) \le 3$.\\
Therefore $\loc^*(S) \subset {\proj^4}^*$ is an irreducible threefold which contains a two-parameter
family of planes. It is easy to prove that in this case $\loc^*(S)$ is a hyperplane of $\pf^*$.\\
It follows that $S$ parametrizes a family of planes in $\pf^*$ contained in a 
hyperplane, and hence, by duality, $S$ parametrizes a two-dimensional family of lines passing through a point
of $\pf$, and it is therefore contained in a $\Sigma$-solid.
\end{proof}

\begin{lemma}\label{sec} 
Let $S$ be a surface in $\G(1,3) \subset \proj^5.$ If $\ord(S) \ge 2$ 
or $\cl(S)\ge 2$, then there exist proper secant lines of $S$ which are contained in $\G(1,3)$.
\end{lemma}

\begin{proof} 
Let $p \in \proj^3$ be a general point. The order of $S$ is the number of lines
parametrized by $S$ which pass through $p$. Hence, if $\ord(S) \ge 2$, there exist at least two lines 
$l_1$, $l_2$ parametrized by $S$ containing $p$.
The pencil of lines generated by $l_1$ and $l_2$ corresponds to a line in $\G(1,3)$
joining the points $L_1$, $L_2 \in S$. Since $p$ is general, the general member of the pencil
is not a line parametrized by $S$, and hence the corresponding secant is not contained in $S$.\par
Let $\pi \subset \proj^3$ be a general plane; the class of $S$ is the number of lines
parametrized by $S$ contained in $\pi$. So if $\cl(S) \ge 2$ there exist $l_1$, $l_2 \subset \pi$, and
the pencil of lines generated by $l_1$ and $l_2$ corresponds to a line in $\G(1,3)$
joining the points $L_1$ and $L_2$. Since $\pi$ is general, the general member of the pencil
is not a line parametrized by $S$, and hence the corresponding secant is not contained in $S$.
\end{proof}

\begin{corollary}\label{corsec}
If $S \subset \G(1,3)$ and $\deg S \ge 3$ then there exist proper secant lines of 
$S$ which are contained in $\G(1,3)$.
\end{corollary}

\begin{proposition}\label{biquadric}
Let $\qu \subset \G(1,4) \subset \proj^9$ be a  two-dimensional smooth quadric such that no proper
secant of $\qu$ is contained in $\G(1,4);$ then $\qu$ is contained in a $\G(1,3)$ and has  
bidegree $(1,1)$. 
In particular, $\qu$ parametrizes the family of lines which lie in a hyperplane 
$H \subset \proj^4$ and meet two skew lines $r$, $s \subset H$. 
\end{proposition}

\begin{proof}
We have $2= \deg (\qu)= \ord(\qu)+\cl(\qu);$ by Lemma \ref{zero} we cannot have $\ord(S)=0$.
If $\ord(S)=2$ then $\cl(S)=0$ and the same Lemma yields that $\qu$ is contained in 
a $\Sigma$-solid, and in this case all the lines in the $\Sigma$-solid meet $\qu$ and are 
contained in $\G(1,4)$.
Therefore $\ord(\qu)=1$ and the statement follows by \cite[Main Theorem]{Ran}.
\end{proof}

\begin{proposition}\label{biscroll}
Let $S \subset \G(1,4)$ be a surface of degree three such that no proper
secant of $S$ is contained in $\G(1,4);$ then the  
bidegree of $S$ is $(2,1)$ and $S$ is not contained in any $\G(1,3)$.
\end{proposition}

\begin{proof}
We have $3= \deg (S)= \ord(S)+\cl(S);$ we cannot have $\ord(S)=0$ by Lemma \ref{zero}.
By the same lemma, if $\ord(S)=3$ then $S$ is contained in a $\Sigma$-solid, 
and in this case all the lines in the $\Sigma$-solid are secant to $S$ and lie in $\G(1,4)$.\\
If $S \subset \G(1,3)$ then $S$ has proper secants contained in $\G(1,3)$ by Lemma \ref{sec}.
Moreover if $\ord(S)=1$ then  $S \subset \G(1,3)$ by \cite[Main Theorem]{Ran}.
\end{proof}

\begin{proposition} \label{twoscrolls}
Let $\scroll \subset \G(1,4)$ be a surface of bidegree $(2,1)$
not contained in a subgrassmannian $\G(1,3)$. Then $\scroll$
parametrizes lines which are contained in a family $F_1$ of planes of a quadric cone 
 ${\mathcal C} \subset \proj^4$ with zero-dimensional vertex
and meet a given line $m$ which lies in a plane $\pi_m \in F_2,$ 
where $F_2$ is the other family of planes of $\mathcal C$.
\end{proposition}

\begin{proof}
Identifying $\G(1,4)$ with the Grassmannian $\G(2,4)$ of planes in the dual space ${\proj^4}^*$, $\scroll$ 
can be viewed as a surface which parametrizes a two-dimensional family of planes in ${\proj^4}^*$.
The duality exchanges order and class, so $\scroll$, as a subvariety of $\G(2,4)$,
has bidegree $(1,2)$.
We apply \cite[Main theorem]{Ran} and we have the following description of $\scroll$:  

\ms
\noindent
let $\beta: \Bl_{M^*}(\pf^*) \to \pf^*$ be the blow-up of $\pf^*$ along a plane 
$M^* \subset {\proj^4}^*$. We can write  $\Bl_{M^*}({\proj^4}^*)=\proj_\pu(\E)$, where
$\E:=\Ol_\pu^3 \oplus \Ol_\pu(1)$;
denote by $p$ the projection $\Bl_{M^*}({\proj^4}^*) \to \pu$.\\
Let $\F$ be a quotient of $\E$ with $\rk(\F)=\deg \F =2$ and denote by 
$p_0:=p_{|\proj(\F)}$. 
$$\xymatrix@C=10pt@R=30pt{
& \Bl_{M^*}({\proj^4}^*) \ar[dr]^p \ar[dl]_\beta  &  \\
{\proj^4}^* &&   \pu }$$
Then
$$\scroll=\scroll(M^*, \F):=\{\pi \in \G(2,4)\ |\ \beta(p_0^{-1}(x)) \subset \pi \subset \beta(p^{-1}(x))
 \,\,\,{\rm for~some} \,\,\, x \in \pu\}.$$
Since $\E$ is nef also $\F$ is, so $\F = \Ol_\pu(a)\oplus \Ol_\pu(b)$ with $a, b \ge 0$ and $a+b=2$.
Therefore two cases can occur:

\begin{itemize}
\item[\rm(i)] $a=1$, $b=1$, i.e., $\proj(\F) \simeq \pu \times \pu$. In this case the tautological bundle $\xi_\E$
restricts to $\F$ as $\Ol(1,1)$, so the image $\beta(\proj(\F)) \subset {\proj^4}^*$ is a 
smooth quadric $\qu$.
The plane $M^*$ contains a line in one ruling of the quadric, and
$\scroll(M^*, \F)$ parametrizes planes in ${\proj^4}^*$ which intersect 
$M^*$ along this line and contain  a line belonging to the other ruling of $\qu$.\\
Passing to the dual we have the claimed description of $\scroll$, where $m$ is the dual line to
the plane $M^*$.

\item[\rm(ii)] $a=0$, $b=2$, i.e., $\proj(\F) \simeq \mathbb F_2$. In this case the tautological bundle $\xi_\E$
restricts to $\F$ as $C_0+2f$, so the image $\beta(\proj(\F)) \subset {\proj^4}^*$ is a quadric cone 
whose vertex is a point $h^* \in M^*$, therefore all the planes parametrized by $\scroll$ pass through $h^*$.\\
It follows that all the lines parametrized by $\scroll \subset \G(1,4)$ are contained in the hyperplane $H$,
dual to $h^*$; in particular, $\scroll$ is contained in $\G_H(1,3)$. This contradicts our hypothesis and thus
exclude this case.
\end{itemize}
\vspace{-1 cm}
\end{proof}

\section{Getting started}\label{start}

\begin{remark}
Let $X$ be a Fano fivefold with Picard number $\rho_X \ge 2$ and index $r_X=2$;
then $X$ has pseudoindex two.
In fact, by \cite{ACO}, the generalized Mukai conjecture 
$$\rho_X(i_X-1) \le \dim X$$
holds for a Fano fivefold, hence we have that $i_X$ cannot be a multiple of $r_X=2$. 
\end{remark}

\begin{lemma} \label{onlyblow}
Let $X$ be a Fano fivefold of index two and $\sigma:X \to X'$ a birational extremal contraction of 
$X$ which contracts a divisor to a surface. Then $\sigma$ is a smooth blow-up.
\end{lemma}

\begin{proof}
Let $R_\sigma$ be the extremal ray in $\cone(X)$ corresponding to $\sigma$.
From the fiber locus inequality we have $l(R_\sigma)=2$, since the general fiber of $\sigma$ is 
two-dimensional. Let $A'$ be a very ample line bundle on $X'$; the line bundle $A = H \otimes \sigma^*A'$ is relatively
ample and $K_X+ 2A = 2\sigma^*A'$ is a supporting divisor for $\sigma$.
We can thus apply \cite[Corollary 5.8.1]{AWcon} to get that $\sigma$ is equidimensional and
the statement then follows from \cite[Theorem 5.2]{AOsperays}.
\end{proof}

\ms
\begin{proposition} 
Let $X$ be a Fano fivevold of index two which is the blow-up
of a smooth variety $X'$ along a smooth center $T;$ then the cone of curves of
$X$ is one among those listed in the following table$,$ where $F$ denotes a fiber type extremal ray$,$ $D_i$
denotes a birational extremal ray whose associated contraction contracts a divisor
to an $i$-dimensional variety and $S$ denotes a ray whose associated contraction is small$:$\\

\begin{center}
\begin{tabular}{c||c|c|c|c||c}
 $\rho_X$ \quad & \quad $R_1$ \quad & \quad $R_2$ \quad & \quad $R_3$ \quad & \quad $R_4$ \quad & \\
\hline\hline
$2$  & $F$ & $D_0$  & & & \rm(a)\\
  & $F$ & $D_2$  & & & \rm(b)\\
	& $D_2$ & $D_2$ &  & & \rm(c)\\
	 & $D_2$ & $S$  & & & \rm(d)\\
\hline
$3$ & $F$ & $F$ & $D_2$  & & \rm(e)\\
    & $F$ & $D_2$ & $D_2$  & & \rm(f)\\
\hline
$4$ & $F$ & $F$ & $F$ & $D_2$ & \rm(g)\\
\hline
\end{tabular}
\end{center}~\par
\end{proposition}

\begin{proof}
The result will follow from the list in \cite[Theorem 1.1]{CO},
once we have proved that $X$ has no contractions of type $D_1$.\\
Let $\sigma: X \to X'$ be the blow-up of $X'$ along $T$, let $E$ be the exceptional divisor
and let $l$ be a line in a fiber of $\sigma$. 
Let $H$ be the fundamental divisor of $X$; from the canonical bundle formula 
$$-2H=K_X = \sigma^*K_{X'} + (\codim T-1) E$$
we know that $-2H \cdot l= (\codim T-1) E \cdot l$, so the codimension of $T$ is odd. 
It follows that either $T$ is a surface in $X'$ or $T$ is a point.
\end{proof}

In this paper we will deal with cases  (b), (e) and (f), since the other cases have already
been classified; in particular:

\begin{itemize}
\item in case (a) $X' \simeq \proj^5$ by \cite[Théorème 1]{BCW}.

\item As noted in the introduction of \cite{CO}, for a Fano fivefold of pseudoindex $2$ possessing a quasi-unsplit locally
unsplit dominating family of rational curves is equivalent to have a fiber type elementary contraction,
so, in cases (c) and (d), we can apply \cite[Theorem 1.2]{CO} and see that  either $X' \simeq \proj^5$ and 
$T$ is

\begin{itemize}
\item[(c1)] a Veronese surface,
\item[(c2)] $\proj_\pu(\Ol(1)\oplus \Ol(2))$ embedded in a 
hyperplane of $\proj^5$ by the tautological bundle (a cubic scroll),
\item[(d1)] a two-dimensional smooth quadric (a section of $\Ol(2)$ in a linear 
$\proj^3 \subset \proj^5$),
\end{itemize}

\ms
\noindent
or $X'$ is a del Pezzo manifold of degree five and $T$ is a plane
of bidegree $(0,1)$. This corresponds to case (c3) which arises as the other extremal contraction 
of case (c2); for a detailed description see \cite[Section 3, Example e1]{CO}.\par

\item In case (g) $X' \simeq \pu \times \pu \times \pt$ and $T \simeq \pu \times \pu \times \{p\}$ by 
\cite[Corollary 5.3]{NO}.
\end{itemize}


\section{Case (b)} \label{rho2}

\subsection{Classification of $X'$} 

We will now prove that if $X$ is as in case $(b)$ then $X'$ is either the projective space of dimension five
or a del Pezzo manifold of degree $\le 5$.\par

\ms
Assume throughout the section that $X$ is a Fano fivefold of index two with $-K_X=2H$ and Mori cone
$\cone(X)=\twospan{R_\vartheta}{R_\sigma}$, where $\vartheta:X \to Y$ is a fiber type contraction
and $\sigma:X \to X'$ is a blow-down with center a smooth surface $S \subset X'$ and 
exceptional divisor $E$. 
By \cite[Theorem 1]{Bo} we know that $X'$ is a smooth Fano variety with $\rho_{X'}=1$ and $i_{X'} \ge 2$; 
moreover by the canonical bundle formula 
$$K_X = \sigma^*K_{X'} + 2E$$
we have that $r_{X'}$ is even. 

\begin{lemma}\label{famv} 
Let $V'$ be a minimal dominating family for $X',$  $V$  a family of 
deformations of the strict transform of a general curve in $V'$ and  
$\V$  the Chow family associated to $V$.
Then $E \cdot V=0,$ the family $\V$ is not quasi-unsplit and $-K_{X'} \cdot V'=4$ or $6$.
\end{lemma}

\begin{proof} 
By \cite[II.3.7]{Kob}, the general curve in $V'$ does not intersect $S$, so $E \cdot V = 0$. 
It follows that
\begin{equation}\label{uno}
-K_X \cdot V = -K_{X'} \cdot V' \le \dim X'+1=6.
\end{equation}
The family $V$ is dominating  and it is not extremal, otherwise
$E$ would be non positive on the whole cone of $X$.
This implies by \cite[Lemma 2.4]{CO} that $X$
is rc$\V$-connected; in particular, since $\rho_X=2$, the family
$\V$ is not quasi-unsplit.
Therefore $-K_{X'} \cdot V'=-K_X \cdot V \ge 4$ so,
recalling that $r_{X'}$ is even, the lemma is proved.
\end{proof}

If the anticanonical degree of the minimal dominating family $V'$ is equal to
$6=\dim X'+1$ then $X' \simeq \proj^5$ by \cite[Theorem 1.1]{Kepn} 
(Note that the assumptions of the quoted result are different, but
the proof actually works in our case since for a very general $x' \in X'$ the pointed family 
$(V')_{x'}$ has the properties 1-3 in  \cite[Theorem 2.1]{Kepn}).\par
\ms
We are thus left with the case $-K_{X'} \cdot V'=4$, which requires some more work.\\
First of all we will analyze the families of rational curves on $X$; as a consequence we will prove
that the exceptional divisor $E$ of the blow-up is a Fano manifold and 
that the fiber type extremal contraction of $X$ restricts to an extremal contraction of $E$ with the same target $Y$.
Using the classification of Fano bundles over a surface, given in \cite{SzW1} and \cite{SzW3}
and completed in Section \ref{fb} of the present paper, we will find a line bundle on 
$Y$ whose pullback to $X$ has degree one on the fibers of the blow-up, and this implies the existence
of a line bundle on $X'$ which has degree one on the rational curves of minimal degree in $X'$.
In this way we will be able to show that $X'$ is a del Pezzo manifold.

\begin{lemma} \label{Ecuts} 
Let  $D$ be an effective divisor of $X;$ then $D$ contains curves whose nu\-me\-rical class is in $R_\sigma$.
\end{lemma}

\begin{proof}
We can assume that $D \not = E$, otherwise the statement is trivial.\\
The image of $D$ via $\sigma$ is an effective divisor in $X'$, hence it is ample since $\rho_{X'}=1$;
therefore $\sigma(D) \cap S \not = \emptyset$ and so $D \cap E \neq \emptyset$.
Let $x$ be a point in $D \cap E$ and let $F_x$ be the fiber of $\sigma$ through $x$;
since $\dim F_x=2$ then $D \cap F_x$ contains a curve in $F_x$.
\end{proof} 

\begin{lemma}\label{nofamily} 
Let $W$ be an unsplit family of rational curves on $X$ such that $\loc(W) \subseteq E;$ 
then $[W] \in R_\sigma$.
\end{lemma}

\begin{proof} 
Let $F$ be a fiber of $\sigma$ such that $F \cap \loc(W) \not= \emptyset$; we have 
$\loc(W)_F \subseteq \loc(W) \subseteq E$.
Assume that $[W] \not \in R_\sigma$; we can apply Lemma \ref{locy} to get
$\dim \loc(W)_F = 4$, so in this case $E = \loc(W)_F = \loc(W)$ and $\cone^X(E)=\twospan{[W]}{R_\sigma}$ 
by Proposition \ref{coneS} (c).
\\
It follows that $E$ contains two independent unsplit dominating families, and it is easy
to prove that their degree with respect to $-K_E$ is equal to three;
we can therefore apply \cite[Theorem 1]{Op} and obtain that $E \simeq \pd \times \pd$.
\\
The effective divisor $E$, being negative on $R_\sigma$, must be positive on $R_\vartheta$, so $E$
dominates $Y$; since $\pd \times \pd$ is a toric variety, by \cite[Theorem 1]{OW} we have that 
$Y \simeq \proj^4$. 
\\
Moreover $\vartheta:X \to \pf$ is a $\pu$-bundle by \cite[Corollary 2.15]{NO};
by \cite[Theorem 1.2]{NO} it must be $X \simeq \proj_{\proj^4}(\Ol \oplus \Ol(a))$ with $a=1$ or
$a=3$, and in these cases $X$ is not a blow-up along a surface, a contradiction.
\end{proof}

\begin{lemma}\label{nodiv}
There does not exist on $X$ any unsplit family of rational curves $W$ 
which satisfies all the following conditions$:$

\item[\rm(i)]  $-K_X \cdot W = 2;$
\item[\rm(ii)]  $[W]$ is not extremal in $\cone(X);$
\item[\rm(iii)] $D_W:=\loc(W)$ has dimension $4;$
\item[\rm(iv)]  $\cone^X(D_W) \subset \langle R_\sigma, [W] \rangle$.

\end{lemma}

\begin{proof} Assume by contradiction that such a family exists. In this case we have
$D_W \cdot R_\sigma \ge 0$ (otherwise we would have $D_W=E$ and $[W] \in R_\sigma$
by Lemma \ref{nofamily}, against
assumption (ii)) and $D_W \cdot R_\vartheta >0$ (otherwise $D_W$ would contain curves in $R_\vartheta$,
against assumption (iv));
this implies that $D_W$ is nef, and that it possibly vanishes only on
$R_\sigma$.
\\
By \cite[Corollary 2.15]{NO} the contraction $\vartheta:X \to Y$ is a $\pu$-bundle, 
i.e., $X=\proj_Y(\E=\vartheta_*H)$; by the classification in \cite[Theorem 1.3]{NO} (note that we are in case $\rho_X=2$) 
this is possible only if $Y$ is a Fano manifold of index one and pseudoindex two or three; in fact
in none of the other cases of \cite[Theorem 1.3]{NO} $X$ is the blow-up of a smooth variety
along a (smooth) surface.

Let $V_Y$ be a family of rational curves on $Y$ with $-K_Y \cdot V_Y=i_Y$ and let
$\nu: \pu \to Y$ be the normalization of a curve in $V_Y$; 
the pull-back $\nu^*\E$ splits as 
$\Ol_{\pu}(1) \oplus \Ol_{\pu}(1)$ in case $i_Y=2$, and as 
$\Ol_{\pu}(1) \oplus \Ol_{\pu}(2)$ in case $i_Y=3$.
We have a commutative diagram
$$
\xymatrix@=35pt{
S:=\proj(\nu^*\E)  \ar[r]^{\bar \nu} \ar[d]_{p} & X \ar[d]^{\vartheta}\\
\pu  \ar[r]_{\nu} & Y }
$$

Let $C \subset S$ be a section corresponding to a surjection $\nu^*\E \to \Ol_{\pu}(1) \to 0$,
and let $V_C$ be the family of deformations of $\bar\nu(C)$; since $H \cdot \bar \nu (C)=\Ol_{\proj(\nu^*\E)}(1) \cdot C=1$ 
the family $V_C$ has anticanonical degree two and is unsplit.

\ms
We claim that the numerical class of $W$ lies in the interior of the cone spanned
by  $[V_C]$ and $R_\vartheta$; this is trivial if $[V_C] \in R_\sigma$, so we can assume that
this is not the case.
\\
The cone of curves of $S$ is generated
by the numerical class of a fiber and the numerical class of $C$, i.e.,
$\cone(S) = \twospan{[C]}{[f]}$.
The morphism $\bar \nu$ induces a map $N_1(S) \to N_1(X)$
which allows us to identify $\cone(S)$ with the subcone of $\cone(X)$ generated
by $[V_C]$ and $R_\vartheta$. The divisor $D_W$ is positive on this subcone, hence the effective divisor
$\Gamma=\bar\nu^*D_W$ is ample on $S$.
It follows that $\Gamma$ lies in the interior of $\cone(S)$, hence $\bar\nu (\Gamma)$,
which is a curve in $D_W$, lies in the interior of the cone generated by $[V_C]$ and $R_\vartheta$.
Therefore also $[W]$ lies in the interior of the cone generated by $[V_C]$ and $R_\vartheta$ 
by assumption (iv), and we can write
$$[W]=a[C_\vartheta] +b [V_C], \qquad \mathrm{with} \quad a,b > 0,$$
where $C_\vartheta$ is a minimal curve in $R_\vartheta$.
\\
Intersecting with $H$ we get $a+b=1$, and intersecting with $-\vartheta^*K_Y$
we have
$$-\vartheta^*K_Y \cdot W=bi_Y < i_Y;$$
therefore if $C_W$ is a curve in $W$ we have
 $-K_Y \cdot \vartheta_*(C_W) < i_Y$, a contradiction.
\end{proof}

\begin{proposition}\label{extrcomp}
Let $V'$ be a minimal dominating family for $X',$  $V$ a family of 
deformations of the strict transform of a curve in $V'$ and $\V$ the Chow family associated to $V$.
Assume that $-K_{X'}\cdot V'=4$.
Then any irreducible component of a reducible cycle in $\V$ which
is not numerically proportional to $V$ is a minimal extremal curve.
\end{proposition}

\begin{proof}
Let $\Gamma=\sum \Gamma_i$ be a reducible cycle in $\V$ with $[\Gamma_1] \not = \lambda [V]$;
since $r_X=2$, $\Gamma$ has exactly two irreducible components.
Denote by $W$ and $\W$ their families of deformations, which have anticanonical
degree two and so are unsplit.
\\
Since by Lemma \ref{famv} $E \cdot V = 0$, we can assume that $E \cdot W < 0$, hence by Lemma 
\ref{nofamily} we have 
that $[W] \in R_\sigma$.

As a consequence, note that if $\Gamma'=\Gamma_1'+\Gamma_2'$ is another reducible cycle in $\V$, then
either $\Gamma_1'$ and $\Gamma_2'$ are numerically proportional to $V$ or, 
denoted by $W'$ and $\W'$ their families of deformations, we can assume that
$[W']=[W]$ and $[\W']=[\W]$.

We claim that $[\W]$ is extremal.\par
\medskip
{\bf Case 1} \quad $V$ is not locally unsplit.\par
\medskip
Let $\{\W^i\}_{i=1, \dots,n}$ be the families of deformations of the irreducible components
of cycles in $\V$ such that $[\W^i]=[\W]$; since $V$ is not locally unsplit,
for some index $i$ the family $\W^i$ is dominating.
We can then apply  \cite[Lemma 2.4]{CO}.\par
\medskip
{\bf Case 2} \quad $V$ is locally unsplit.\par
\medskip
Assume by contradiction that $[\W]$ is not extremal.
By the argument in the proof of Case 1 we have that $\W^i$ is not dominating for every $i$.
By inequality \ref{iowifam} (a) we have that $\dim \loc(\W^i)=3$ or $4$; we distinguish two cases:

\begin{itemize}
\item[(i)] There exists an index $i$ such that $\dim \loc(\W^i)=4$.

\ms
Let $D=\loc(\W^i)$; if $D \cdot V=0$ then $D$ is negative on an extremal ray of 
$\cone(X)$, hence on $R_\sigma$, but this implies $D=E$,
against Lemma \ref{nofamily}.
\\
Therefore $D \cdot V  > 0$, hence $D \cap \loc(V_x) \neq \emptyset$ for a general $x \in X$.  
Since we are assuming that $V$ is locally unsplit we have that $\dim \loc(V_x) \ge 3$
and $\cone^X(\loc(V_x))=\onespan{V}$ by Proposition \ref{coneS} (b), 
so $\dim\loc(\W^i)_{\loc(V_x)} \ge 4$ by Lemma \ref{locy} (b) and $D =\loc(\W^i)_{\loc(V_x)}$.
It follows by \cite[Lemma 1]{Op} that every curve in $D$ 
can be written as $a C_V + b C_{\W^i}$ with $a \ge 0$, $C_V$ a curve contained
in $\loc(V_x)$ and $C_{\W^i}$ a curve
in $\W^i$. 
Therefore $\cone^X(D) \subset \twospan{R_\sigma}{[\W^i]}$, but this is excluded by Lemma
\ref{nodiv}.

\item[(ii)] For every $i$ we have $\dim \loc(\W^i)=3$.

\ms
By inequality \ref{iowifam} (a) we have $\dim \loc(\overline W_x) =3$ for every $x \in \loc(\W)$.
\\
Let $\Omega = \cup_i(\loc(W^i) \cup \loc(\W^i))= E \cup_i \loc(\W^i)$, and take a point
$y$ outside $\Omega$; since $X$ is rc$\V$-connected we can join $y$ and $\Omega$ 
with a chain of cycles in $\V$.
Let $C$ be the first irreducible component of these cycles which meets $\Omega$.
\\
Clearly $C$ cannot belong to any family $W^i$ or $\W^i$ because it is not contained in $\Omega$,
so it belongs either to $V$ or to a family $\lambda V$ which is numerically
proportional to $V$; by \cite[Lemma 9.1]{ACO} we have 
that either $C \subset \loc(V_z)$ for some $z$ such that $V_z$ is unsplit or
$C \subset \loc(\lambda V)$.
\\
Moreover, since $E \cdot V =0$ the intersection
$C \cap \Omega$ is contained in $\Omega \setminus E$.
Let $t$ be a point in $C \cap \Omega$ and let $\Omega_j =\loc(\W^j)$ be the irreducible 
component of $\Omega$ which contains $t$.
If $C \subset \loc(V_z)$ we have $\dim (\loc(V_z) \cap \Omega_j) \ge 1$, against the fact that 
$N_1^X(V_z)=\onespan{[V]}$ and $N_1^X(\Omega_j)=\onespan{[\W^j]}$.
\\
If else $C \subset \loc(\lambda V)$ we have that $\dim \loc(\lambda V)_{\Omega_j} \ge 4$ by Lemma \ref{locy} (b)
and that $\cone^X(\loc(\lambda V)_{\Omega_j}) \subset \langle [\lambda V], R_\vartheta \rangle$
by \cite[Lemma 1]{Op}; this is clearly impossible if $\loc(\lambda V)_{\Omega_j}=X$,
and it contradicts Lemma \ref{Ecuts} if $\dim \loc(\lambda V)_{\Omega_j}=4$. 
\end{itemize}
Finally, since $-K_X \cdot W^i =-K_X \cdot \W^i=2$ we also have that the curves of
$W^i$ and $\W^i$ are minimal in $R_\sigma$ and $R_\vartheta$ respectively.
\end{proof}

\begin{corollary}\label{Vinmezzo}
In the assumptions of Proposition \ref{extrcomp}$,$ denoting as usual by $C_\sigma$ and $C_\vartheta$ minimal
rational curves in the rays $R_\sigma$ and $R_\vartheta$, we have, in $\cone(X)$, $[V]=[C_\sigma]+[C_\vartheta]$;
in particular we have $H \cdot C_\vartheta=1$.
\end{corollary}

\begin{proposition}\label{X'dp}
Let $V'$ be a minimal dominating family for $X',$ let $V$ be a family of 
deformations of the strict transform of a curve in $V'$ and assume that $-K_{X'}\cdot V'=4$.
Then $E$ is a Fano manifold and  $X'$ is a del Pezzo manifold.
\end{proposition}

\begin{proof} 
By Lemma \ref{famv} we have $E \cdot V=0$, hence 
$E \cdot C_\vartheta= -E \cdot C_\sigma=1$ by Corollary \ref{Vinmezzo};
It follows that
$$(-K_X  -E)\cdot C_\sigma = 2+1 = 3$$
$$(-K_X  -E)\cdot C_\vartheta = 2-1 = 1,$$
hence $-K_X -E$ is ample on $X$ by Kleiman criterion.
By adjunction $-K_E=(-K_X-E)_{|E}$ is ample on $E$
and $E$ is a Fano manifold. \par
\ms
We note that $E$ contains curves of $R_\vartheta$: otherwise the 
fiber type contraction $\vartheta$ would be a $\pu$-bundle by \cite[Lemma 2.13]{NO}, and since
$E \cdot C_\vartheta = 1$ it follows that $E$ would be  a section of
$\vartheta$, against the fact that $\rho_Y=1$ and $\rho_E = \rho_S+1 \ge 2$.
\\
Consider the divisor $D=H-E$: it is nef and vanishes on $R_\vartheta$, so it is a supporting divisor
for $\vartheta$. The restriction $D_{|E}$ is nef but not ample, since $E$ contains curves of $R_\vartheta$, 
so $D_{|E}$ is associated to an extremal face of $\cone(E)$ and to an extremal contraction 
$\vartheta_E:E \to Z$ and  we have a commutative diagram:
$$\xymatrix@=40pt{
E \ar[r] \ar[dr]|{\vartheta_{|E}} \ar[d]_{\vartheta_E}  & X \ar[d]^\vartheta \\
Z \ar[r] & Y }$$\par

\ms
We will prove that, for every $m \in \mathbb N$ the restriction map $H^0(X,mD) \to H^0(E,mD_{|E})$ is an isomorphism,
hence $\vartheta_{|E}=\vartheta_E$ and $Z=Y$.
Consider the exact sequence
$$\xymatrix{0 \ar[r] & \Ol_X(mD-E) \ar[r] & \Ol_X(mD) \ar[r] & \Ol_E(mD_{|E}) \ar[r] & 0}.$$
Since $E$ is not contracted by $\vartheta$ we have that $h^0(mD-E)=0$; 
moreover, we can write 
$$mD-E=K_X+(m-1)D+3H-2E.$$ 
By Kleiman criterion $3H-2E$ is ample on $X$ and, being $(m-1)D$ nef, the divisor $(m-1)D+3H-2E$ is ample, too. 
By the Kodaira Vanishing Theorem $h^1(mD-E)=0$.
\\
We have proved that $E$ is a Fano manifold, and we know that it has a $\pd$-bundle structure
over $S$, i.e., $E \simeq \proj_S(\E)$ with $\E$ a Fano bundle of rank three over $S$. 
This implies that $S$ is a del Pezzo surface.\par

\ms
Let $L_Y$ be the ample generator of $\pic(Y)$; by Proposition \ref{fbodps}, Proposition \ref{fanoquadr}
and the classification in \cite{SzW1}, the pull-back of
$L_Y$ has degree one on the fibers of the  $\pd$-bundle.
\\
The line bundle $H-E$ has degree two on the fibers of the $\pd$-bundle and is trivial on the fibers of 
$\vartheta$, hence $H-E=2\vartheta^*L_Y$ and so $H-\vartheta^*L_Y$ is trivial on the fibers of $\sigma$,
i.e., $H-\vartheta^*L_Y=\sigma^*H_{X'}$ for some $H_{X'} \in \pic(X')$.
By the canonical bundle formula we have
\begin{equation}\label{index4}
-\sigma^*K_{X'}=-K_X+2E=2(H+E)=4H -4 \vartheta^*L_Y=4\sigma^*H_{X'},
\end{equation}
i.e., $r_{X'}=4$ and so $X'$ 
is a del Pezzo fivefold.
\end{proof}

\begin{corollary}
By the classification of del Pezzo manifolds given by Fujita \cite{Fuj}$,$ denoting 
by $d:=H_{X'}^5$ the degree of $X'$ and recalling that $\rho_{X'}=1,$ 
we have the following possibilities$:$

\item[\rm(i)] if $d=1$ then $X' \simeq V_1$ is a degree six hypersurface in the weighted projective 
space $\proj(3,2,1,\ldots,1);$
\item[\rm(ii)] if $d=2$ then $X' \simeq V_2$ is a double cover of $\proj^5$  branched along a smooth 
quartic hypersurface$;$
\item[\rm(iii)] if $d=3$ then $X' \simeq V_3$ is a cubic hypersurface in $\proj^{6};$
\item[\rm(iv)] if $d=4$ then $X' \simeq V_4$ is the complete intersection of two quadrics in $\proj^{7};$
\item[\rm(v)] if $d=5$ then $X' \simeq V_5$ is a linear section of the grassmannian $\G(1,4) \subset \proj^9$.

\end{corollary}


\subsection{Classification of $S$} 

\begin{theorem}\label{p5} 
If $X' \simeq \proj^5$ then $S$ is as in Theorem \ref{main}, cases \rm{(b1)-(b6)}.
\end{theorem}

\begin{proof}
Let $H$ be a hyperplane of $\proj^5$,  let $\tl H \subset X$ be its strict transform via $\sigma$ and
let $\ac = \sigma^*H$. 
We know that $\tl H$ is an effective divisor different from $E$, hence it is nef; moreover if $S \subset H$
we can write $\tl H= \ac-kE$ with $k >0$.
\\
Let $\Gamma$ be a proper bisecant of $S$, and let $\tl \Gamma$ be its strict transform; if $S \subset H$ 
we have
$$0 \le \tl H \cdot \tl \Gamma \le 1-2k;$$
it follows that $S$ has no proper bisecants, i.e., $S$ is a linear subspace of $\proj^5$ and
we are in case (b1).
If else $S$ is not contained in any hyperplane, note that $S$ cannot be the Veronese surface,
since the blow-up of $\proj^5$ along a Veronese surface has two birational contractions; therefore
the secant variety of $S$ fills $\proj^5$.\par

\ms
Let $l$ be a line in $\proj^5$ not contained in $S$ and $\tilde l$ its strict transform; we have
$$-K_X \cdot \tilde l = \sigma^*\Ol_{\proj^5}(6) \cdot \tilde l - 2E\cdot \tilde l = 6 - 2(\sharp(S \cap l))>0;$$
therefore if $l$ is a proper bisecant of $S$ we have $-K_X \cdot \tilde l = 2$; moreover 
$S$ cannot have (proper) trisecant lines.
\\
In the notation of \cite{Bau}, the condition on the trisecants is equivalent to the fact 
that the trisecant variety of $S$ (which consists of all lines contained in $S$ and of the
proper trisecants) is contained in $S$, so by
the description in \cite{Bau} (see in particular Theorem 7, Section 4 and Appendix A2)
we have the possibilities (b2)-(b6).\par

\medskip
We now show that in all these cases the blow-up of $X'$ along $S$ is a Fano manifold
with the prescribed cone of curves.
\\  
The linear system $\mathcal L=|\Ol_{\proj^5}(2) \otimes {\cal I}_S|$
of the quadrics in $\proj^5$ containing $S$ has $S$ as its base locus scheme
(see \cite{Io}), so
$\sigma^*\mathcal L$ defines a morphism $\vartheta:X \to \proj(\mathcal L)$.
\\
Since $2 \ac - E$ is nef and vanishes on the strict transforms
of the bisecants of $S$ it follows that the numerical class of these 
curves is extremal in $\cone(X)$, and since $-K_X$ is positive on these curves we can conclude that
$X$ is a Fano manifold.
\\
Moreover since $S$ is neither degenerate nor the Veronese surface, the bisecants to $S$ cover 
$\proj^5$ and so $\vartheta$ is of fiber type.
\end{proof}

\begin{lemma}\label{condnum} 
Assume that $X'$ is a del Pezzo fivefold.
Let  $H_{X'} = \Ol_{X'}(1)$ and $H_S = (H_{X'})_{|S}$. Then
\begin{enumerate}
\item[\rm(i)] If $\dim Y=2$ then $H_S^2= \deg X' =-K_S\cdot H_S$.
\item[\rm(ii)] If $\dim Y=3$ then $\deg X'= -K_S \cdot H_S$ and $\deg X'-H_S^2 \ge 2$.
\item[\rm(iii)] If $\dim Y=4$ then $\deg X' > - K_S \cdot H_S$.
\end{enumerate}
\end{lemma}

\begin{proof}
Denote by $\Norm$ the normal bundle of $S$ in $X'$ and by $\Norm^*$ the conormal bundle; 
let $C = \det \Norm^* \in \pic(S)$.
Recall that $E = \proj_S(\Norm^*)$ and that $-E_{|E}=\xi_{\Norm^*}$. 
\\
Let $\ac=\sigma^*H_{X'}$; we have
$$\ac^5 = (H_{X'})^5=\deg X'=:d,$$
and since the intersection of three or more sections of a very ample multiple of $H_{X'}$ 
does not meet $S$ we have also
$$\ac^4E=\ac^3E^2=0.$$
Then we have
\begin{equation*}
\begin{array}{l}
K_S = (K_{X'}+\det \Norm)_{|S}=-4H_S - C,\\
\ac^2E^3= (\ac^2E^2)_{|E}=H_S^2,\\
\ac E^4=(\ac E^3)_{|E} =(-\ac\xi_{\Norm^*}^3)_{|E}= -C \cdot H_S,
\end{array}
\end{equation*}
Let $L:=\ac -E$; from the above equalities it follows that
\begin{equation}\label{for}
L^4\ac  = \ac^5-4\ac ^2E^3+\ac E^4= d +K_S \cdot H_S;
\end{equation}
\begin{equation}\label{for2}
L^3\ac^2  = \ac^5 -\ac^2E^3 = d - H_S^2.
\end{equation}
By Corollary \ref{Vinmezzo} we have that $H \cdot C_\vartheta =1$; then equation (\ref{index4}) yields that 
$\ac \cdot R_\vartheta = E \cdot R_\vartheta=1$, hence $L$ is trivial on 
the fibers of $\vartheta$ and therefore $L=\vartheta^*L_Y$.

\begin{enumerate}
\item[(i)]
If $\dim Y = 2$ we have  $L^4\ac = L^3\ac^2 = 0$, so it follows from (\ref{for}) and (\ref{for2}) that
$$0 = d +K_S \cdot H_S = d -H_S^2.$$

\item[(ii)]
If  $\dim Y = 3$ then  $L^4\ac=0$, and so by \ref{for} we have
$$ d +K_S \cdot H_S =0.$$
The contraction $\vartheta$ is a quadric fibration (see Definition \ref{defcontrazioni})
and $\ac_{|F}=\Ol_F(1)$ for a general fiber 
$F$ of $\vartheta$; hence $L^3 \ac^2 = (L_Y^3)(\ac_F^2) \ge 2$, and \ref{for2} yields that
$$d - H_S^2 \ge 2.$$

\item[(iii)]
Finally, if $\dim Y = 4$ the general fiber $F$ of $\vartheta$ is one-dimensional and $\ac \cdot F=1$,
hence $L^4 \ac = L_Y^4 > 0$; again by \ref{for} we have that
$$d + K_S \cdot H_S >0.$$
\end{enumerate}
\vspace{-0.5cm}
\end{proof}

\begin{lemma}
If $\dim Y >2$ then $S$ is $\pd,$ a smooth quadric $\qu$ or the ruled surface $\mathbb F_1,$ i.e.
the blow-up of $\pd$ at a point.
\end{lemma}

\begin{proof} 
By Proposition  \ref{X'dp} $E$ is a Fano manifold and, by the proof of the same Proposition, 
we know that the restriction $\vartheta_{|E}:E \to Y$
is an extremal contraction of $E$. Moreover, by the classification in Proposition \ref{fbodps} 
we know that for every del Pezzo surface $S_k$ with $k \ge 2$
the exceptional divisor $E$ is isomorphic to
$S_k \times \pd$, and in this case $E$ has no maps on a variety with Picard number 
one and dimension greater than two.  
\end{proof}

\begin{theorem} \label{dpclass}
If $X'$ is a del Pezzo fivefold then the pairs $(X',S)$ are  as in Theorem \ref{main}, cases \rm{(b7)-(b13)}.
\end{theorem}

\begin{proof}
The contraction $\vartheta:X \to Y$ is supported by $\ac - E$, and is the
resolution of the rational map $\theta:X' \xymatrix{ \ar@{-->}[r] &}Y$ defined by the linear system 
$\mathcal L :=\sigma_*|\vartheta^*L_Y|$, where $L_Y$ is the ample generator
of $\pic(Y)$; since $|\vartheta^*L_Y|$ is base point free
we have $\Bs \mathcal L \subseteq S$; on the other hand
$\mathcal L \subseteq |H_{X'} \otimes {\cal I}_S|$, therefore $\Bs \mathcal L \supseteq S$
and so $\Bs \mathcal L = S$.
\\
It follows that the strict transforms of curves of degree one
with respect to $H_{X'}$ which meet $S$ are contracted by $\vartheta$. 
Moreover, since $\ac - E$ is nef, no curves of degree one with respect 
to $H_{X'}$ and not contained in $S$ can meet $S$ in more than one point.

\begin{itemize}
\item If $\dim Y = 2$ then $\vartheta$ is equidimensional and by \cite[Corollary 1.4]{AWcon} we have that 
$Y$ is smooth; moreover $\rho_Y=1$ and $Y$ is dominated by a Fano manifold, so $Y \simeq \proj^2$.
\\
Therefore $\dim \mathcal L =3$, so $S$ is the complete intersection
of three general sections in $|H_{X'}|$ and we are in case (b7).

\item In case  $\dim Y = 3$, if $S\simeq \pd$ then $H_{S}\simeq \Ol_\pd(a)$, with $a >0$. 
By Lemma \ref{condnum} (ii) we have $d=-K_{\pd} \cdot H_{\pd}=3a$; recalling that $d \le 5$
we find $H_S = \Ol_{\pd}(1)$ and $d =3$ (case (b8)).
\\
If $S \simeq \pu \times \pu$ then $H_{S}\simeq \Ol_{\pu \times \pu}(a,b)$, with $a,b >0$.
By Lemma \ref{condnum} (ii) we have $d=-K_{\pu \times \pu} \cdot H_{\pd}=2a +2b$;
recalling that $d \le 5$ we find $H_S = \Ol_{\pu \times \pu}(1,1)$ and $d =4$ (case (b10)).
\\
For $S \simeq \mathbb F_1$ we have $-K_{\mathbb F_1} \cdot C \ge 5$ 
for every ample $C \in \pic(\mathbb F_1)$, equality holding if and only if $C=C_0+2f$;
hence, by Lemma \ref{condnum} (ii) we have  $d=-K_{\mathbb F_1} \cdot H_{\mathbb F_1}=5$
and $H_S=C_0+2f$.
Since all the bisecants of $S$ which are contained in $\G(1,4)$ are also contained in a linear section $V_5$,
it follows by Proposition \ref{biscroll} that $S$ is as in case (b13).

\item Finally, in case $\dim Y = 4$ we can apply Lemma \ref{condnum} (iii) and get: if 
$S \simeq \pd$ then $H_S=\Ol(1)$ and $H_S^2=1$, so $d=4$ (case (b9)) or $d=5$; in the latter case,
being $\vartheta$ of fiber type, we exclude the case of a plane of bidegree $(0,1)$ in view of Remark
\ref{blowupaplane} and we are in case (b11).
\\
If $S \simeq \pu \times \pu$ the bound $-K_S \cdot H_S \ge 4$ gives $H_S=\Ol(1,1)$ and $H_S^2=2$, 
hence $d=5$; in this case $S$  has bidegree $(1,1)$  by Proposition \ref{biquadric}  and we are in
case (b12). 
The center of the blow-up cannot be $\mathbb F_1$ since 
 $-K_{\mathbb F_1} \cdot H_{\mathbb F_1} \ge 5$, which contradicts Lemma \ref{condnum} (iii).
\end{itemize}
 
\ms
We show now that in all these cases the blow-up of $X'$ along $S$ is a Fano manifold
with the prescribed cone of curves.
\\ 
Let $(X',S)$ be a pair as in the theorem
and denote by $H_{X'}$ the fundamental divisor of $X'$.
We claim that the linear system $|H_{X'} \otimes {\cal I}_S|$  has $S$ as its base locus scheme;
this is clear apart from cases (b10), which is described in Proposition \ref{quadr}, and (b12) and (b13), 
which are treated in Proposition \ref{appe}.
\\
Therefore the linear system $|\sigma^*H_{X'} - E|$ defines a morphism 
$\vartheta:X \to \proj(|\sigma^*H_{X'} - E|)$.
Since $\sigma^*H_{X'} - E$ is nef and vanishes on the strict transforms
of the rational curves of degree one in $X'$ which meet $S$, it follows that the numerical class of these 
curves is extremal in $\cone(X)$.
Being $-K_X$ positive on these curves, we can conclude that $X$ is a Fano manifold.
Finally, since the curves of degree one with respect to $H_{X'}$ which meet $S$ cover $X'$  
we have that $\vartheta$ is a fiber type contraction.
\end{proof}

\begin{proposition} \label{quadr}
Let ${\qu}$ be a smooth two-dimensional quadric in 
$V_4 \subset \proj^7$. Then $\qu$ is the intersection of $V_4$ and the hyperplanes of $\proj^7$
which contain $\qu$.
\end{proposition}

\begin{proof}
Let ${\qu}$ be a smooth two-dimensional quadric in 
$V_4 = \quadr\, \cap\, \quadr' \subset \proj^7$, and let 
$\Lambda_{\qu}^3$ be the three-dimensional linear subspace of $\proj^7$
which contains ${\qu}$.
\\
We claim that $\Lambda_{\qu}^3$ is contained in one of the two quadrics $\quadr$,
$\quadr'$.
From \cite[Proposition 2.1]{Reidth} we know that the intersection of two quadrics  
is smooth if and only if there exist coordinates in $\proj^n$ such that
$$\quadr = \left\{ \sum x_i^2=0 \right\}\qquad \quadr' = \left\{ \sum \lambda_ix_i^2=0 \right\}$$
with $\lambda_i \neq \lambda_j$ for every $i \neq j$. So assume by contradiction that 
$\Lambda_{\qu}^3 \not \subset \quadr \cup \quadr'$; 
in this case $\Lambda_{\qu}^3 \cap \quadr = \Lambda_{\qu}^3 \cap \quadr' 
= {\qu}$, so it must be 
$$\left(\sum (1-\lambda_i)x_i^2\right)_{|\Lambda_{\qu}^3} \equiv 0.$$
But there is at most one index $i$ such that $\lambda_i=1$, so the kernel of the quadratic form 
$\sum (1-\lambda_i)x_i^2$ is at most one-dimensional and we reach a contradiction.
\end{proof}

\begin{proposition} \label{appe}
Let $S$ be a smooth two-dimensional quadric of bidegree $(1,1)$ or a surface of bidegree $(2,1)$ 
not contained in a $\G(1,3)$, in $V_5 \subset \proj^8$. Then $S$ is the intersection of $V_5$ 
and the hyperplanes of $\proj^8$ which contain $S$.
\end{proposition}

\begin{proof} Since $V_5$ is an hyperplane section of $\G(1,4)$ we will show that 
$S \subset \G(1,4) \subset \proj^9$ is the intersection of $\G(1,4)$ 
and the hyperplanes of $\proj^9$ which contain $S$, by finding explicitly its equations.\\
By Proposition \ref{biquadric}, if $S$ is a quadric of bidegree $(1,1)$, then
it parametrizes lines in $\pf$ which meet two given skew lines $r$, $s$.
Up to a change of coordinates in $\proj^4$, we can assume that
$r$ and $s$ have equations
$$r: \ x_0=x_1=x_2=0, \quad \quad s: \  x_0=x_3=x_4=0,$$
so $H$ is the hyperplane of equation $x_0=0$; in this case
 the equations of $S$ in $\G$ are 
$$\left\{\begin{array}{l}y_0 = \ldots = y_4 = y_9 = 0 \\
y_5y_8 =  y_6y_7  
\end{array}\right.$$
and $S$ is the intersection of $\G$ with the three-dimensional linear subspace 
$\Lambda^3 \subset \proj^9$ of equations 
$$y_0= \ldots = y_4 = y_9 = 0.$$
Let now $S \subset \G$ be a surface of bidegree $(2,1)$ not contained in a $\G(1,3)$, 
as described in Proposition \ref{twoscrolls}. Up to a coordinate change in $\proj^4$, assume that
${\mathcal C}$ is the cone of vertex $(0:0:0:0:1)$ on the quadric of equations 
$$x_0x_2=x_1x_3, \qquad x_4=0,$$
and that $m$ is the line of equations $x_0=x_1=x_4=0$.
\noindent
The two families of planes contained in ${\mathcal C}$ have equations
$$F_1 = 
\begin{cases} 
\lambda x_0 = \mu x_1 \\ 
\lambda x_3 = \mu x_2 
\end{cases} \qquad \qquad
F_2 = 
\begin{cases} \lambda x_0 = \mu x_3 \\ 
\lambda x_1 = \mu x_2 
\end{cases},$$
and $m$ lies in the plane $\pi_m \in F_2$ of equations $x_0=x_1=0$.
\noindent
The equations of the scroll $S \subset \G$ are
$$\left\{\begin{array}{l}
y_0 =y_3=y_6= y_7 =0 \\
y_1=y_5 \\
y_1^2=y_2y_4 \\
y_1y_8 =  y_4y_9\\
y_1 y_9 = y_2 y_8 \\
\end{array} \right.$$
\noindent
In particular, $S$ is the intersection of $\G$ with the four-dimensional linear
space $\Lambda_S^4$ of equations $y_0=y_3=y_6=y_7=0$, 
$y_1=y_5$.
\end{proof}

\section{Cases (e)-(f)}\label{rho3}

{\bf Setup.} \quad Throughout the section, let $X$ be a Fano fivefold whose cone of curves 
is as in cases (e)-(f),
and let $\sigma:X \to X'$ be an extremal contraction of $X$ which is the blow-up
of $X'$ along a smooth surface.

\begin{proposition}\label{samelocus}
Let $X$ be as above.
Then either $X=\proj_{\pd \times \pd}(\Ol \oplus \Ol(1,1))$
or $X'$ is a Fano manifold of even index.
\end{proposition}

\begin{proof}
Let $E$ be the exceptional locus of $\sigma$; by \cite[Proposition 3.4]{Wicon} $X'$
is a Fano manifold unless $E$ contains the exceptional locus of another
extremal ray; this is clearly possible only if $X$ has another birational contraction, i.e., in case (f).
Note that in this case both the birational contractions of $X$
are smooth blow-ups by Lemma \ref{onlyblow}.
\\
Let $\os$ be the other blow-up contraction of $X$, 
denote by $R_\sigma$ and $R_{\os}$ the extremal rays corresponding
to $\sigma$ and $\os$ and by $R_\vartheta$ the extremal ray
corresponding to the fiber type contraction $\vartheta:X \to Y$.
\\
Let $F$ be a fiber of $\sigma$; by Lemma \ref{locy} (a) we have 
$\dim \loc(R_\os)_F \ge 4$, hence
$E=\loc(R_\os)_F$ and $\cone^X(E)=\langle R_\sigma,R_\os \rangle$ 
by Proposition \ref{coneS}.
Moreover $E \cdot R_\sigma < 0$ and $E \cdot R_\os < 0$, hence $E \cdot R_\vartheta >0$  
and $\vartheta$ is a $\pu$-bundle by \cite[Corollary 2.15]{NO}.
We can thus apply \cite[Theorem 1.1]{NO}, noting that the only Fano manifold in the list
given in that result with two birational contractions with the same exceptional
locus is $X=\proj_{\pd \times \pd}(\Ol \oplus \Ol(1,1))$.
\\
The claim about the index of $X'$ follows from the canonical bundle
formula for $\sigma$.
\end{proof}

\begin{lemma}\label{Epos} 
Let $X$ be a Fano fivefold whose cone of curves is as in case $(f);$
denote by $R_\sigma$ and $R_\os$ the divisorial extremal rays of $\cone(X),$ by
$R_\vartheta$ the fiber type extremal ray and by $E$ (resp. $\overline E$) the exceptional locus of 
$R_\sigma$ (resp. $R_\os$).
Then either $E \cdot R_\vartheta >0,$ or $\overline E \cdot R_\vartheta>0,$.
\end{lemma}

\begin{proof}
Consider a minimal horizontal dominating family $V$ for $\vartheta$.\par

\ms
{\bf Claim.}\quad 
\emph{The numerical class of $V$ belongs to a two-dimensional extremal face of $\cone(X)$
which contains $R_\vartheta$.}

\medskip
If $V$ is unsplit, since $\rho_X=3$ the claim follows from \cite[Lemma 2.4]{CO}.

\medskip
Denote by $V_\vartheta$ the family of 
deformations of a minimal curve in $R_\vartheta$.
If $V$ is not unsplit, for a general $x \in \loc(V)$ we have that
$\dim \loc(V_x) \ge 3$ by Proposition \ref{iowifam},  $\cone^X(\loc(V_x)) = \onespan{V}$ 
by Proposition \ref{coneS} and $\dim \loc(V_\vartheta,V)_x \ge 4$ by Lemma \ref{locy} (c).
Call $D = \loc(V_\vartheta,V)_x$: then $\cycl^X(D)=\twospan{R_\vartheta}{V}$ by  \cite[Lemma 1]{Op}, 
so $D$ is a divisor since $\rho_X = 3$.
\\
It cannot be $D \cdot R_\vartheta > 0$, otherwise
we could write $X=\cloc(V_\vartheta,V)_x$ 
and we would have $\rho_X=2$; so it must be $D \cdot R_\vartheta = 0$.
This implies that $D$ is positive on a birational ray, say $R_\sigma$, 
hence $\dim (D \cap F) \ge 1$ for every fiber $F$ of $\sigma$; 
since $\cycl^X(D)=\twospan{R_\vartheta}{V}$ 
and $\cone^X(F) = \onespan{R_\sigma}$, the claim is proved.

\medskip
It follows that $E \cdot R_\vartheta>0$: in fact, if $E \cdot R_\vartheta=0$
then $E\cdot V <0$, since curves of $V$ are not contracted by $\vartheta$ and so 
they do not belong to $R_\vartheta$.
But then we would have $\loc(V) \subset E$ and $V$ would not be
dominating for $\vartheta$, a contradiction.
\end{proof}

\begin{proposition}\label{difflocus}
Let $X$ be a Fano fivefold whose cone of curves is as in cases $(e)-(f),$
and let $\sigma:X \to X'$ be the blow-up
of $X'$ along a smooth surface$;$ assume that $E$ is positive
on a fiber type extremal ray of $X$.
\\
If $X'$ is a Fano manifold$,$ then either $X' \simeq \pu \times \quadr^4,$ and in this case 
either $S \simeq \pu \times l$
with $l$ a line in $\quadr^4$ or $S \simeq \pu \times \Gamma$ with $\Gamma$ a conic not
contained in a plane $\pi \subset \quadr^4,$  
or $X'$ is a $\proj^3$-bundle over $\pd$ and $S$ dominates $\pd$ via the bundle projection.
\end{proposition}

\begin{proof}
Let $R_\vartheta$ be the extremal ray on which $E$ is positive,
and let $\vartheta:X \to Y$ be its associated contraction; let 
$\psi:X \to W$ be the contraction of the face 
spanned by $R_\sigma$ and $R_\vartheta$. Then 
$\psi$ factors through $\sigma$ and a 
morphism $\theta:X' \to W$, and we have a commutative diagram
$$\xymatrix@=40pt{
X \ar[d]_\sigma \ar[dr]^{\psi} \ar[r]^\vartheta & Y \ar[d]\\
X'  \ar[r]_{\theta} & W}$$

The contractions $\sigma$ and $\psi$ have connected fibers, so the same is true
for $\theta$; moreover $W$ is a normal variety with $\rho_W=\rho_{X'}-1$ and
$\dim W < \dim X'$.
\\
It follows that $\theta$ is an extremal elementary fiber type contraction of the 
Fano manifold $X'$; denote by $R_\theta$ the corresponding extremal ray in $\cone(X')$.

Let $V'_\theta$ be a dominating family of rational curves whose numerical class
belongs to $R_\theta$ and whose degree with respect to some ample line bundle 
is minimal among the degrees of the families with this
property. In particular, by the minimality assumption, such a family is locally unsplit.
Let $V$ be the family of deformations of
the strict transform in $X$ of a general curve in $V'_\theta$.
\\
Since curves of $V$ are contracted by $\psi$,
the numerical class of $V$ in $\cone(X)$ lies in the face spanned
by $R_\sigma$ and $R_\vartheta$.
\\
By \cite[II.3.7]{Kob}, the general curve in $V'_\theta$ does not intersect the center 
$S$ of the blow-up, so $E \cdot V = 0$; it follows that $[V] \not \in R_\vartheta$.
Clearly we cannot have $[V]  \in R_\sigma$, being $E \cdot R_\sigma <0$,
so the class $[V]$ does not generate an extremal ray of $X$.
In particular, since $V$ is dominating and $X$ has no small
contractions, $V$ cannot be unsplit in view of \cite[Lemma 2.29]{CO}, hence
$$4 \le -K_X \cdot V=-K_{X'} \cdot V'_\theta.$$
For a general $x \in X'$ we have, by Proposition \ref{iowifam} (b) that
$\dim \loc(V'_\theta)_x \ge 3$, so a general fiber of  $\theta$ is at least 
three-dimensional and $\dim W \le 2$.

\ms
If $\dim W=1$ then the contraction of the extremal ray of $X$ different
from $R_\sigma$ and $R_\vartheta$ is a $\pu$-bundle by \cite[Corollary 2.15]{NO}
(take a fiber of $\psi$ for $D$).
Now we apply \cite[Lemma 4.1]{NO}, to get that $X$ is a product with $\pu$
as a factor; looking at the classification table in \cite[Appendix]{NO} we find that
the only products with $\rho_X=3$ and a blow-down contraction of type $D_2$ are
$X \simeq \pu \times \Bl_l(\quadr^4)$ or
$X \simeq \pu \times \Bl_\Gamma(\quadr^4)$; the
description of $X'$ and $S$ follows.

\ms
If $\dim W=2$ we claim that $X'$ is a $\proj^3$-bundle over $\pd$. We would like to use
\cite[Lemma 2.18]{NO}, but we do not know that the length of the ray $R_\theta$ is $\dim X'-1$.
However we notice that, in the proof of the quoted result, the assumption on the length
is used only to prove that the general fiber of the contraction is a projective space, so we will prove
in a different way that this is the case in our situation.\\
Let $x$ be a general point in $X'$ and denote by $F_x$ the fiber of $\theta$ containing $x$; 
by Proposition \ref{iowifam} (b) we have $\dim \loc(V'_\theta)_x \ge 3$, hence $F_x = \loc(V'_\theta)_x$.
Moreover, since $V'_\theta$ is locally unsplit, by Proposition \ref{coneS} (b), 
we have $\rho_{F_x}=1$. Now we can conclude $F_x \simeq \proj^3$ either by the classification
of Fano threefolds or by applying \cite[Theorem 1.1]{Kepn} as in the proof of Lemma \ref{famv}.\\ 
Therefore, by the proof of \cite[Lemma 2.18]{NO}, $X'$ is a $\proj^3$-bundle over $\pd$; 
$E$ is positive on the fiber type ray $R_\vartheta$, so the image via $\sigma$ 
of every curve in $R_\vartheta$ is a curve contracted by $\theta$ which meets $S$.
Since $\vartheta$ is a fiber type contraction, we know that curves in $R_\vartheta$
dominate $X$, hence curves contracted by $\theta$ which meet $S$ 
dominate $X'$. Therefore $S$ dominates $\pd$.
\end{proof}

\begin{theorem} 
Let $X$ be a Fano fivefold whose cone of curves is as in cases \rm{(e)-(f)},
and let $\sigma: X \to X'$ be the blow-up of $X'$ along a smooth surface $S$. 
Then the pairs $(X',S)$ are as in Theorem \ref{main}, cases \rm{(e1)-(e4)} or \rm{(f1)-(f4)}.
\end{theorem}

\begin{proof}
By Proposition \ref{samelocus}, either $X \simeq \proj_{\pd \times \pd}(\Ol \oplus \Ol(1,1))$
and therefore $(X',S)$ is as in case (f1) or  we can apply 
Proposition \ref{difflocus}: in fact, in case (e) the positivity of $E$ 
on a fiber type ray of $\cone(X)$ is trivial, otherwise it follows from Lemma \ref{Epos}.
\\
Therefore either $(X',S)$ is as in cases (e1)-(e2) or, 
up to exchange $\sigma$ with $\os$, we have that
$X'$ is a $\proj^3$-bundle over $\pd$.
\\
In this case, the classification in \cite{SzW1} yields that $X'$ is either the blow-up of $\proj^5$ 
along a plane $\pi_1$ or $X' \simeq \proj_{\pd}(T\pd(-1) \oplus \Ol^{\oplus 2})$.
Considering the exact sequence
$$\shse{\Ol_\pd(-1)}{\Ol_\pd^{\oplus 5}}{T\pd(-1) \oplus \Ol_\pd^{\oplus 2}}$$
we see that $X'=  \proj_{\pd}(T\pd(-1) \oplus \Ol^{\oplus 2})$ embeds in $\pd \times \proj^4$
as a section of $\Ol(1,1)$.\par
\ms
Let $l \subset X'$ be a line in a fiber of the $\proj^3$-bundle not contained in $S$,
and let $\tilde l \subset X$ be its strict transform; by the canonical bundle formula
\begin{equation}\label{lines}
-K_X \cdot \tilde l = -\sigma^*K_{X'} \cdot \tilde l -2E \cdot \tilde l \le 4-2 \#(S \cap l);
\end{equation}
since $X$ is Fano it must be $\#(S \cap l) \le 1$.
\\
Let $R_\oth \subset \cone(X')$ be the extremal ray of $X'$ not associated to the $\proj^3$-bundle 
contraction. Let $C$ be a minimal extremal curve in $R_\oth$ not contained in $S$
and let $\tl C \subset X$ be its strict transform.
Again by the canonical bundle formula
$$-K_X \cdot \tl C = -\sigma^*K_{X'} \cdot \tl C -2E \cdot \tl C \le 2-2 \#(S \cap C),$$
hence $S \cap C  = \emptyset$.
Therefore, if $S$ meets a two-dimensional fiber $F_\oth$ of $\oth$ then $S=F_\oth$.

\begin{itemize}
\item
In case $X' \simeq \Bl_{\pi_1}(\proj^5)$, the map $\oth$ is the blow-up map, so 
denoted by $E'$ the exceptional divisor of $\oth$ we have that either $S$ is a fiber of $\oth$
and we are in case (f2), or $S \cap E'= \emptyset$; in particular $S$ cannot meet a fiber of the 
$\proj^3$-bundle in a curve.
\\
In the first case, $X$ has another blow-down contraction 
$\os : X \to \Bl_p(\proj^5)$, whose center is the strict transform of a plane passing through $p$;
this corresponds to case (f3).
In fact, $X$ can be described as follows: let
$Y$ be the blow-up of $\pf$ along a line, let $E_Y$ be the exceptional divisor, 
let $H_Y$ be the pullback of $\Ol_\pf(1)$ and
let $\E= (2H_Y +E_Y) \oplus (3H_Y +E_Y)$. Then $X= \proj_{Y}(\E)$, and the following diagram
shows the extremal contractions of $X$:
$$
\xymatrix@=20pt{
\pd  & & \Bl_{\pi_1}(\proj^5) \ar[ll]_\theta \ar[rd]^{\oth} &\\
\Bl_l(\pf) \ar[d] \ar[u] &  X \ar[l]^(.35)\vartheta 
 \ar[ru]^\sigma \ar[rd]_{\os} & & \proj^5\\
 \pf  & &  \Bl_{p}(\proj^5)  \ar[ll] \ar[ru]  & }
$$
In case $S \cap E'= \emptyset$, 
equation (\ref{lines}) yields that $S$ is a section of the $\proj^3$-bundle contraction of $X'$; 
therefore it corresponds to a surjection $\Ol^3 \oplus \Ol(1) \to \Ol(1)$,
the image of $S$ in $\proj^5$ is a plane $\pi_2$ not meeting $\pi_1$ and we are in case (f4).
In this case $X \simeq \proj_{\pd \times \pd}(\Ol(0,1) \oplus \Ol(1,0))$.

\item
If $X' \simeq \proj_{\pd}(T\pd(-1) \oplus \Ol^{\oplus 2})$ the contraction $\oth$ is of fiber type;
it follows that $S$ is the union of all the fibers
of $\oth$ which have nonempty intersection with $S$ itself.
\\
In particular, either $S$ is a two-dimensional fiber of $\oth$, i.e.,a 
section corresponding to a surjection $T\pd(-1) \oplus \Ol^{\oplus 2} \to \Ol$,
and we are in case (e3), or $\oth$ is a $\pu$-bundle and $S$ contains a one-parameter 
family of fibers  isomorphic to $\pu$.
\\
In this last case, the restriction of $\oth$ to $S$ is a morphism from $S$ to a curve, and therefore
$S \not \simeq \pd$; so $S$ cannot be a section of the natural projection $p:X' \to \pd$.
By equation (\ref{lines}) the restriction of $p$ to $S$ is a birational
morphism $p_{|S}:S \to \pd$, and the only surface which is birational to $\pd$ and has a 
morphism on a curve all whose fibers are $\simeq \pu$ is the Hirzebruch surface $\mathbb F_1$. 
In particular, the exceptional curve of $S$ is a line in a fiber of $p$, therefore
$\oth(\mathbb F_1)=\oth(C_0)$ is a line $l \subset \proj^4$ and
$S$ is the intersection of the pullback of three hyperplanes in $\pf$
meeting along $l$ (case (e4)).

\ms
To conclude, we prove the effectiveness of $X$ in these last two cases:
in case (e3) let $Y$ be a general member of $|\Ol_{\pd \times \pt}(1,1)|$, and let
$\E= \Ol_Y(1,1) \oplus \Ol_Y(1,2)$; then $X \simeq \proj_Y(\E)$, as proved in \cite[Proposition 7.3]{NO}, 
and  $X$ is a $\pu$-ruled Fano manifold.
\\
In case (e4) $X$ can be realized as follows: let $Z=\Bl_l(\proj^4)$, and let $H_Z$
be the pullback of $\Ol_\pf(1)$; then $X$ is a general section
in the linear system $|p_1^*\Ol_\pd(1)+p_2^*H|$ in $\pd \times Z$, where $p_1$ and $p_2$ 
denote the projections  onto the factors. 
\end{itemize}
\vspace{-0.87 cm}
\end{proof}

\subsection*{Acknowledgements}

We would like to thank the referees for pointing out inaccuracies and mistakes in the
first version of the paper, as well as for their comments, which were of great help
in improving the exposition.


\address{Elena Chierici\\
Dipartimento di Matematica \\
via Sommarive 14, I-38050 Povo (TN) \\
Italy}
{e.chierici@email.it}
\address{Gianluca Occhetta\\
Dipartimento di Matematica \\
via Sommarive 14, I-38050 Povo (TN) \\
Italy}
{gianluca.occhetta@unitn.it}

\end{document}